\documentclass[twoside]{amsart}
\usepackage[centertags]{amsmath}

\usepackage{amsfonts}
\usepackage{amssymb}
\usepackage{amsthm}
\usepackage{eucal}
\usepackage[cmtip,all]{xy}
\usepackage{graphicx}

\def\I{^{\mathrm I}}
\def\II{^{\mathrm{I\!I}}}

\newdir{ >}{{}*!/-20pt/@{>}}
\newdir{> }{{}*!/10pt/@{>}}
\newdir{>> }{{}*!/10pt/@{>>}}

\newtheorem{thm}[equation]{Theorem}
\newtheorem{cor}[equation]{Corollary}
\newtheorem{lem}[equation]{Lemma}

\theoremstyle{definition}
\newtheorem{defn}[equation]{Definition}
\theoremstyle{remark}
\newtheorem{rem}[equation]{Remark}
\newtheorem{exm}[equation]{Example}
\numberwithin{equation}{section}

\newcommand{\abs}[1]{\left\vert#1\right\vert}
\newcommand{\set}[1]{\left\{#1\right\}}
\newcommand{\ina}[1]{\left\{#1\right\}}

\newcommand{\eps}{\varepsilon}
\newcommand{\To}{\longrightarrow}

\newcommand{\A}{\mathcal{A}}

\newcommand{\C}[1]{\mathbf{#1}}
\newcommand{\D}[1]{\mathcal{#1}}

\newcommand{\EZDIAG}[5]{\xymatrix@C+=2.5cm{*+[r]{#1}
\ar@(u,l)_(0.62){\displaystyle #5}[]
\ar@<.5ex>^-{#3}[r]&\ar@<.5ex>^-{#4}[l]#2}}

\def\ni{{nil}}
\def\abb{{ab}}

\def\r{\rightarrow} 
\def\rr{\Rightarrow} 

\def\ner{\operatorname{Ner}}

\def\ho{\operatorname{Ho}}

\def\diag{\operatorname{Diag}}

\def\st{\stackrel} 

\def\ul{\underline} 

\def\ca{\circledcirc}
\def\co{\bar{\otimes}}

\def\coker{\operatorname{Coker}}
\renewcommand{\ker}{\operatorname{Ker}}

\def\Z{\mathbb{Z}}

\def\L{\Omega}

\newcommand{\grupo}[1]{\langle #1\rangle}

\newcommand{\vc}{\Box}    

\begin{document}

\title{The $1$-type of a Waldhausen $K$-theory spectrum}%
\author{Fernando Muro \and Andrew Tonks}%
\address{Universitat de Barcelona, Departament d'$\grave{\text{A}}$lgebra i Geometria, Gran Via de les Corts Catalanes, 585, 08007 Barcelona, Spain}
\email{fmuro@ub.edu}
\address{London Metropolitan University, 166--220 Holloway Road, London N7 8DB, UK}
\email{a.tonks@londonmet.ac.uk}

\thanks{The first author was partially supported
by the project MTM2004-01865 and the MEC postdoctoral fellowship
EX2004-0616, and the second by the MEC-FEDER grant MTM2004-03629.}

\begin{abstract}
We give a small functorial algebraic model for the 2-stage Postnikov
section of the $K$-theory spectrum of a Waldhausen category and use
our presentation to describe the multiplicative structure with
respect to biexact functors.
\end{abstract} \maketitle

\section*{Introduction}

Waldhausen's $K$-theory of a category $\C{C}$ with cofibrations and
weak equivalences \cite{akts} extends the classical notions of
$K$-theory, such as the $K$-theory of rings, additive categories and
exact categories.

In this paper we give an algebraic model $\D{D}_*\C{C}$ for the 1-type
$P_1K\C{C}$ of the $K$-theory spectrum $K\C{C}$. This model consists
of a diagram of groups
$$\xymatrix{&
(\D{D}_0\C{C})^\abb\otimes
(\D{D}_0\C{C})^\abb\ar[d]_{\grupo{\cdot,\cdot}}&&
\\K_1\C{C}\;\ar@{^{(}->}[r]&
\D{D}_1\C{C}\ar[r]^\partial&\D{D}_0\C{C}\ar@{->>}[r]\ar@{-->}@/_20pt/[lu]_-H&K_0\C{C}.}$$
in which the bottom row is exact.

The important features of our model are the following:
\begin{itemize}
\item It is {\em small}, as it has generators given by the objects,
weak equivalences and cofiber sequences of the category $\C{C}$.
\item It has {\em minimal nilpotency degree}, since both groups
$\D{D}_0\C{C}$ and $\D{D}_1\C{C}$ have nilpotency class 2.
\item It encodes the 1-type in a {\em functorial }way,
and the homotopy classes of morphisms $\D{D}_*\C{C}\r\D{D}_*\C{D}$
and $P_1K\C{C}\r P_1K\C{D}$ are in bijection.
\end{itemize}

{}From this structure we can recover the homomorphism $\cdot\eta\colon
K_0\C{C}\otimes\Z/2\r K_1\C{C}$, which gives the action of the Hopf
map in the stable homotopy groups of spheres, in the following way
$$
a\cdot \eta\;\;=\;\; \grupo{a,a},\quad a\in K_0\C{C}.
$$

The extra structure given by the quadratic map $H$ is used to
describe the behaviour of $\D{D}_*$ with respect to biexact functors
between Waldhausen categories $\C C\times \C D\r \C E$. In
particular the classical homomorphisms
\begin{eqnarray*}
K_0\C{C}\otimes
K_0\C{D}&\To&K_0\C{E},\\
K_0\C{C}\otimes K_1\C{D}&\To&K_1\C{E},\\
K_1\C{C}\otimes K_0\C{D}&\To&K_1\C{E},
\end{eqnarray*}
may be obtained from our model $\D D_*$.

The object $\D{D}_*\C{C}$ is a {\em stable quadratic module} in the
sense of \cite{ch4c}. This object is defined below by a presentation in terms of generators and relations
in the spirit of Nenashev, who gave a model for $K_1$ of an exact
category in \cite{K1gr}. A stable quadratic module is a particular
case of a strict symmetric categorical group, or more generally of a commutative monoid in the category of crossed complexes,
which were first introduced by Whitehead in \cite{chII}. The
monoidal structure for crossed complexes was defined in \cite{tpcc}.

To obtain our presentation of $\D{D}_*$
we introduce the {\em total crossed complex} $\Pi X$ of a
bisimplicial set $X$, and show that there is an Eilenberg--Zilber--Cartier
equivalence
\[\EZDIAG{\pi\diag(X)}{\Pi(X)}{}{}{}\]
generalizing \cite[Section 2]{hnaf} and \cite{EZcc}.
This is then applied to the bisimplicial set given by the nerve of
Waldhausen's $wS.$ construction~\cite{akttsI}.
An explicit step-by-step translation from the total complex $\Pi X$ to our model
$\D{D}_*$ is presented.

\section{The algebraic $1$-type $\D{D}_*\C{C}$ of the $K$-theory spectrum $K\C{C}$}\label{1}

We begin by defining the algebraic structure which the model
$\D{D}_*\C{C}$ will have.

\begin{defn}\label{ob}
A \emph{stable quadratic module} $C_*$ is a diagram of group
homomorphisms
$$C_0^\abb\otimes C_0^\abb\st{\omega}\To C_1\st{\partial}\To C_0$$
such that given $c_i,d_i\in C_i$, $i=0,1$,
\begin{itemize}
\item $\partial\omega(\ina{c_0}\otimes\ina{d_0})=[d_0,c_0]$,
\item $\omega(\ina{\partial (c_1)}\otimes\ina{\partial (d_1)})=[d_1,c_1]$,
\item $\omega(\ina{c_0}\otimes\ina{d_0}+\ina{d_0}\otimes\ina{c_0})=0$.
\end{itemize}
Here $[x,y]=-x-y+x+y$ is the commutator of two elements $x,y\in K$
in any group $K$, $K^\abb$ is the abelianization of $K$, and
$\ina{x}\in K^\abb$ is the element represented by $x\in K$. In
order to simplify the notation we will use the following
convention throughout the whole paper
\begin{eqnarray*}
\grupo{c_0,d_0}&=&\omega(\set{c_0}\otimes\set{d_0}).
\end{eqnarray*}
We will also write $\partial_{C_*}$, $\omega_{C_*}$ and
$\grupo{\cdot,\cdot}_{C_*}$ for the structure maps of $C_*$ if we
want to distinguish between different stable quadratic modules.
\end{defn}

Stable quadratic modules were introduced in \cite[Definition
IV.C.1]{ch4c}. Notice, however, that we adopt the opposite
convention for the homomorphism $\omega$.

\begin{rem}\label{basics}
The laws of a stable quadratic module imply that $C_0$ and $C_1$
are groups of nilpotency class $2$. Indeed, given $x,y,z\in C_0$
we have
\begin{equation*}
[x,[y,z]]\quad=\quad\partial\omega(\ina{[y,z]}\otimes\ina{x})\quad=\quad
\partial\omega(0\otimes\ina{x})\quad=\quad0.
\end{equation*}
Moreover, given $a,b,c\in C_1$
\begin{equation*}
[a,[b,c]]\,=\,\omega(\ina{\partial([b,c])}\otimes\ina{\partial(a)})\,=\,
\omega(\ina{[\partial(b),\partial(c)]}\otimes\ina{\partial(a)})\,=\,
\omega(0\otimes\ina{\partial(a)})\,=\,0.
\end{equation*}
The homomorphism $\omega$ is central since
\begin{equation*}
\begin{array}{r}
[a,\omega(\ina{y}\otimes\ina{z})]\;=\;
\omega(\ina{\partial\omega(\ina{y}\otimes\ina{z})}\otimes\ina{\partial(a)})\;=\;
\omega(\ina{[z,y]}\otimes\ina{\partial(a)})\\=\;\omega(0\otimes\ina{\partial(a)})\;=\;0.
\end{array}
\end{equation*}
Furthermore, $\ker\partial\subset C_1$ is a central subgroup since
given $d\in C_1$ with $\partial(d)=0$
\begin{equation*}
[a,d]\;=\;\omega(\ina{\partial(d)}\otimes\ina{\partial(a)})\;=\omega(0\otimes\ina{\partial(a)})\;=\;0.
\end{equation*}
\end{rem}

As usual, one can define stable quadratic modules in terms of
generators and relations. In this way we define the algebraic
model $\D{D}_*\C{C}$ in Definition \ref{LA} below. Free stable
quadratic modules, and also stable quadratic modules defined by
generators and relations in degrees zero and one, can be
characterized up to isomorphism by obvious universal properties.
In Appendix \ref{append} we give explicit constructions of the
groups and the structure homomorphisms of a stable quadratic
module defined by generators and relations.

\begin{exm}
We give an easy example of a presentation of a stable quadratic module which
serves to highlight the difference with group presentations, and also the suitability
for stable homotopy theory. Let $C_*$ be the stable quadratic module with just one generator $\gamma$, in degree zero. In this case $C_0$ is clearly the infinite cyclic group generated by $\gamma$, whereas $C_1$ is isomorphic to $\Z/2$, generated by the element $\grupo{\gamma,\gamma}$.
\end{exm}

We assume the reader has certain familiarity with Waldhausen categories and related concepts. We refer to \cite{wiak} for the basics, see also \cite{akts}.

\begin{defn}\label{LA}
Let $\C{C}$ be a Waldhausen category with distinguished zero object
$*$. Cofibrations and weak equivalences are denoted by
$\rightarrowtail$ and $\st{\sim}\rightarrow$, respectively. A
generic cofiber sequence is denoted by
$$A\rightarrowtail B \twoheadrightarrow B/A.$$

We define $\D{D}_*\C{C}$ as the stable quadratic module generated in
dimension zero by the symbols
\begin{itemize}
\item $[A]$ for any object in $\C{C}$,
\end{itemize}
and in dimension one by
\begin{itemize}
\item $[A\st{\sim}\r A']$ for any weak equivalence,
\item $[A\rightarrowtail B\twoheadrightarrow B/A]$ for any cofiber
sequence,
\end{itemize}
such that the following relations hold.
\begin{enumerate}
\item $\partial[A\st{\sim}\r A']=-[A']+[A]$.
\item $\partial[A\rightarrowtail B\twoheadrightarrow B/A]=-[B]+[B/A]+[A]$.
\item $[*]=0$.
\item $[A\st{1_A}\r A]=0$.
\item $[A\st{1_A}\r A \twoheadrightarrow *]=0$, $[*\rightarrowtail A\st{1_A}\r
A]=0$.
\item For any pair of composable weak equivalences $ A\st{\sim}\r B\st{\sim}\r C$,
$$[A\st{\sim}\r C]=[B\st{\sim}\r C]+[A\st{\sim}\r B].$$
\item For any commutative diagram in $\C{C}$ as follows
$$\xymatrix{A\;\ar@{>->}[r]\ar[d]^\sim&B\ar@{->>}[r]\ar[d]^\sim&B/A\ar[d]^\sim\\A'\;\ar@{>->}[r]&B'\ar@{->>}[r]&B'/A'}$$
we have
\begin{eqnarray*}
[A\st{\sim}\r A']+[B/A\st{\sim}\r B'/A']&&\\+\grupo{[A],-[B'/A']+[B/A]}&=&
-[A'\rightarrowtail
B'\twoheadrightarrow B'/A']\\&&+[B\st{\sim}\r B']\\&&+[A\rightarrowtail B\twoheadrightarrow B/A].
\end{eqnarray*}
\item For any commutative diagram consisting of four obvious cofiber sequences in $\C{C}$ as follows
$$\xymatrix{&&C/B\\&B/A\;\ar@{>->}[r]&C/A\ar@{->>}[u]\\A\;\ar@{>->}[r]&B\;\ar@{>->}[r]\ar@{->>}[u]&C\ar@{->>}[u]}$$
we have
\begin{eqnarray*}
[B\rightarrowtail C\twoheadrightarrow C/B]&&\\+[A\rightarrowtail
B\twoheadrightarrow B/A]&=&[A\rightarrowtail C\twoheadrightarrow
C/A]\\
&&+[B/A\rightarrowtail C/A\twoheadrightarrow C/B]\\
&&+\grupo{[A],-[C/A]+[C/B]+[B/A]}.
\end{eqnarray*}
\item For any pair of objects $A, B$ in $\C{C}$
$$\grupo{[A],[B]}=
{}-[B\st{i_2}\rightarrowtail A\vee B\st{p_1}\twoheadrightarrow A]
+[A\st{i_1}\rightarrowtail A\vee B\st{p_2}\twoheadrightarrow B]
.$$
Here $$\xymatrix{A\ar@<.5ex>[r]^-{i_1}&A\vee
B\ar@<.5ex>[l]^-{p_1}\ar@<-.5ex>[r]_-{p_2}&B\ar@<-.5ex>[l]_-{i_2}}$$
are the natural inclusions and projections of a 
coproduct
in $\C{C}$.
\end{enumerate}
\end{defn}

\begin{rem}
Notice that relations (4) and (7) imply that if (9) holds for a
particular choice of the coproduct $A\vee B$
then it holds for any
other choice $A\vee' B$, since any two coproducts are canonically
isomorphic by an isomorphism which preserves the natural
inclusions and projections of the factors.
\end{rem}

\begin{defn}\label{mor}
A \emph{morphism} $f\colon C_*\r D_*$ in the category $\C{squad}$
of stable quadratic modules is given by group homomorphisms
$f_i\colon C_i\r D_i$, $i=0,1$, such that given $c_i,d_i\in C_i$
\begin{itemize}
\item $\partial_{D_*} f_1(c_1)=f_0\partial_{C_*}(c_1)$,
\item $\grupo{f_0(c_0),f_0(d_0)}_{D_*}=f_1\grupo{c_0,d_0}_{C_*}$.
\end{itemize}
The \emph{homotopy groups} of $C_*$ are
$$\pi_1C_*=\ker\partial\text{  and  }\pi_0C_*=\coker\partial.$$ A
\emph{weak equivalence} in $\C{squad}$ is a morphism which induces
isomorphisms in homotopy groups. The \emph{homotopy category}
$$\ho\C{squad}$$ is obtained from $\C{squad}$ by inverting weak
equivalences.
\end{defn}

Let $\C{WCat}$ be the category of Waldhausen categories as above and
exact functors. 
The construction $\D{D}_*$ defines a functor
$$\D{D}_*\colon\C{WCat}\To\C{squad}.$$
For an exact functor $F\colon \C{C}\r\C{D}$ the stable quadratic module morphism
$\D{D}_*F\colon\D{D}_*\C{C}\r\D{D}_*\C{D}$
is given on generators by
\begin{eqnarray*}
(\D{D}_0F)([A])&=&[F(A)],\\
(\D{D}_1F)([A\st{\sim}\r A'])&=&[F(A)\st{\sim}\r F(A')],\\
(\D{D}_1F)([A\rightarrowtail B\twoheadrightarrow B/A])&=&
[F(A)\rightarrowtail F(B)\twoheadrightarrow F(B/A)].
\end{eqnarray*}

Let $\ho\C{Spec}_0$ be the homotopy category of connective spectra.
In Lemma \ref{equi} below we define a functor
$$\lambda_0\colon\ho\C{Spec}_0\To\ho\C{squad}$$
together with natural isomorphisms
$$\pi_i\lambda_0X\cong\pi_iX,\;\;i=0,1.$$ This functor induces an equivalence of
categories
$$\lambda_0\colon\ho\C{Spec}_0^1\st{\sim}\To\ho\C{squad},$$
where $\ho\C{Spec}_0^1$ is the homotopy category of spectra with
trivial homotopy groups in dimensions other than $0$ and $1$.

The naive algebraic model for the $1$-type of the algebraic $K$-theory spectrum $K\C{C}$ of a Waldhausen category $\C{C}$ would be $\lambda_0K\C{C}$. However this stable quadratic module is much bigger than $\D{D}_*\C{C}$ and it is not directly defined in terms of the basic structure of the Waldhausen category $\C{C}$. This makes meaningful the following theorem, which is the main result of this paper.

\begin{thm}\label{main}
Let $\C{C}$ be a Waldhausen category. There is a natural isomorphism
in $\ho\C{squad}$
$$\D{D}_*\C{C}\st{\cong}\To\lambda_0K\C{C}.$$
\end{thm}

This theorem shows that the model $\D{D}_*\C{C}$ satisfies the
functoriality properties claimed in the introduction. It also shows
that the exact sequence of the introduction is available. The theorem
will be proved in section four.

{}From a local point of view the $1$-type of a connective spectrum is determined up to non-natural isomorphism by the first $k$-invariant. We now establish the connection between this $k$-invariant and the algebraic model $\D{D}_*\C{C}$.

\begin{defn}
The \emph{$k$-invariant} of a stable quadratic module $C_*$ is the homomorphism
$$k\colon\pi_0C_*\otimes\Z/2\To\pi_1C_*,\;\;k(x\otimes 1)=\grupo{x,x}.$$
\end{defn}

Given a connective spectrum $X$ the $k$-invariant of $\lambda_0X$ coincides with the action of the Hopf map
$0\neq\eta\in\pi_1^s\cong\Z/2$ in the stable stem of the sphere.
$$\xymatrix{\pi_0X\otimes\pi_1^s\ar[r]\ar[d]_\cong&\pi_1X\ar[d]^\cong\\
\pi_0\lambda_0X\otimes\Z/2\ar[r]^-k&\pi_1\lambda_0X}$$
See Lemma \ref{equi} below. We recall that the action of $\eta$ coincides with the first Postnikov invariant of $X$. This is
used to derive the following corollary of Theorem \ref{main}.

\begin{cor}
The first Postnikov invariant of the spectrum $K\C{C}$
$$\cdot\eta\colon K_0\C{C}\otimes\Z/2\To K_1\C{C}$$ is defined by $$[A]\cdot\eta=[\tau_{A,A}\colon A\vee A\st{\cong}\r A\vee
A],$$ where $\tau_{A,A}$ is the automorphism which exchanges the factors of a coproduct $A\vee A$ in $\C{C}$.
\end{cor}

\begin{proof}
Consider the following commutative diagram
$$\xymatrix{A\;\ar@{>->}[r]^-{i_2}\ar@{=}[d]&A\vee A\ar@{->>}[r]^-{p_1}\ar[d]^{\tau_{A,A}}_\cong&A\ar@{=}[d]\\
A\;\ar@{>->}[r]^-{i_1}&A\vee A\ar@{->>}[r]^-{p_2}&A}$$ where $i_j$
and $p_j$ are the natural inclusions and projections of the two
factors of the coproduct $A\vee A$, $j=1,2$. By relations (4) and
(7) in Definition \ref{LA}
\begin{eqnarray*}
0&=& -[A\st{i_1}\rightarrowtail A\vee A\st{p_2}\twoheadrightarrow
A]+[\tau_{A,A}\colon A\vee A\st{\cong}\r A\vee
A]+[A\st{i_2}\rightarrowtail A\vee A\st{p_1}\twoheadrightarrow A].
\end{eqnarray*}
By relation (1) in Definition \ref{LA} we have
$\partial[\tau_{A,A}]=0$, so $[\tau_{A,A}]$ is central in
$\D{D}_1\C{C}$ by Remark \ref{basics}, therefore
\begin{eqnarray*}
0&=& [\tau_{A,A}\colon A\vee A\st{\cong}\r A\vee
A]-[A\st{i_1}\rightarrowtail A\vee A\st{p_2}\twoheadrightarrow
A]+[A\st{i_2}\rightarrowtail A\vee A\st{p_1}\twoheadrightarrow
A]\\
&=&[\tau_{A,A}\colon A\vee A\st{\cong}\r A\vee A]-\grupo{[A],[A]}.
\end{eqnarray*}
Here we use relation (9) in Definition \ref{LA} for the last
equality. This proves the corollary.
\end{proof}

The next corollary can be easily obtained from the previous one by using again the relations defining $\D{D}_*$ and matrix arguments as for
example in the proof of
\cite[Proposition 2.1 (iv)]{attI}.

\begin{cor}
Let $\C{A}$ be a Waldhausen category which is additive. Then the first Postnikov invariant of the spectrum $K\C{A}$
$$\cdot\eta\colon K_0\C{A}\otimes\Z/2\To K_1\C{A}$$ is defined by $$[A]\cdot\eta=[-1_A\colon A\st{\cong}\r A].$$
\end{cor}

\section{The multiplicative properties of $\D{D}_*$}

In order to describe the multiplicative properties of $\D{D}_*\C{C}$
with respect to biexact functors we would need a symmetric monoidal
structure on $\C{squad}$ which models the smash product of spectra.
Unfortunately such a monoidal structure does not exist and we need to enrich $\D{D}_*\C{C}$ with an
extra structure map $H$,
$$\xymatrix{
(\D{D}_0\C{C})^\abb\otimes
(\D{D}_0\C{C})^\abb\ar[d]_{\grupo{\cdot,\cdot}}&
\\
\D{D}_1\C{C}\ar[r]^\partial&\D{D}_0\C{C},\ar@/_20pt/[lu]_-H}$$ so
that the diagrams
\begin{eqnarray*}
\D{D}_{0}^{sg}\C{C}&=&\left(\xymatrix{\D{D}_0\C{C}\ar@<.5ex>[r]^-{H}&
\ar@<.5ex>[l]^-{\grupo{\partial,\partial}}(\D{D}_0\C{C})^\abb\otimes
(\D{D}_0\C{C})^\abb}\right),
\\
\D{D}_{1}^{sg}\C{C}&=&\left(\xymatrix{\D{D}_1\C{C}\ar@<.5ex>[r]^-{H\partial}&
\ar@<.5ex>[l]^-{\grupo{\cdot,\cdot}}(\D{D}_0\C{C})^\abb\otimes
(\D{D}_0\C{C})^\abb}\right),
\end{eqnarray*}
are square groups in the sense of
\cite{ecg}.

\begin{defn}
A \emph{square group} $M$ is a diagram
$$M_e\mathop{\rightleftarrows}\limits^{H}_P M_{ee}$$
where $M_e$ is a group, $M_{ee}$ is an abelian group, $P$ is a
homomorphism, $H$ is a quadratic map, i.e.\  the symbol
$$(x|y)_H=H(x+y)-H(y)-H(x),\;\;x,y\in M_e,$$
is bilinear, and the following identities hold, $a\in M_{ee}$,
\begin{eqnarray*}
(P(a)|x)_H&=&0,\\
(x|P(a))_H&=&0,\\
P(x|y)_H&=&[x,y],\\
PHP(a)&=&P(a)+P(a).
\end{eqnarray*}

Note that $(\cdot|\cdot)_H$ induces a homomorphism
\begin{equation}\label{qw}(\cdot|\cdot)_H\colon\coker P\otimes\coker
  P\To M_{ee}.
\end{equation}
Moreover
$$T=HP-1\colon M_{ee}\To M_{ee}$$
is an involution, i.e.\ a homomorphism with $T^2=1$, and
$$\Delta\colon\coker P\To X_{ee}\colon x\mapsto (x|x)_H-H(x)+TH(x)$$
defines a homomorphism.

A \emph{morphism} $f\colon M\r N$ in the category
of square
groups is given by group homomorphisms $f_e\colon M_e\r N_e$,
$f_{ee}\colon M_{ee}\r N_{ee}$ commuting with $H$ and $P$.

A \emph{quadratic pair module} $f\colon M\r N$ is a square group
morphism such that $M_{ee}=N_{ee}$ and $f_{ee}$ is the identity.

Morphisms in the category $\C{qpm}$ of quadratic pair modules are
defined again by homomorphisms commuting with all operators.
\end{defn}

A stable quadratic module $C_*$ is termed \emph{$0$-free} if
$C_0=\grupo{E}^\ni$ is a free group of nilpotency class $2$, see
the appendix. Here $E$ is the basis.

\begin{lem}
Let $C_*$ be a $0$-free stable quadratic module with
$C_0=\grupo{E}^\ni$ and let $H\colon C_0\r\Z[E]\otimes\Z[E]$ be
the unique quadratic map such that $H(e)=0$ for any $e\in E$ and
$(x|y)_H=y\otimes x$ for $x,y\in C_0$. Then
\begin{eqnarray*}
C_0^{sg}&=&\left(\xymatrix{C_0\ar@<.5ex>[r]^-{H}&
\ar@<.5ex>[l]^-{\grupo{\partial,\partial}}C_0^\abb\otimes
C_0^\abb}\right),
\\
C_1^{sg}&=&\left(\xymatrix{C_1\ar@<.5ex>[r]^-{H\partial}&
\ar@<.5ex>[l]^-{\omega}C_0^\abb\otimes C_0^\abb}\right),
\end{eqnarray*}
are square groups. Moreover, the homomorphism $\partial\colon C_1\r C_0$ defines a
quadratic pair module
$$C_1^{sg}\To C_0^{sg}.$$
\end{lem}

The square group $C_0^{sg}$ in this lemma will also be denoted
by $\Z_\ni[E]$ or just $\Z_\ni$ if $E$ is a singleton, as in \cite{qaI}.

The stable quadratic module $\D{D}_*\C{C}$ defined
in the previous section is $0$-free. The basis of $\D{D}_0\C{C}$ is
the set of objects in $\C{C}$, excluding the zero object $*$.
In particular $\D{D}_{0}^{sg}\C{C}$ and $\D{D}_{1}^{sg}\C{C}$ above
are square groups and $\D{D}_*\C{C}$ is endowed with the structure of a
quadratic pair module. Moreover, the morphisms induced by exact
functors are compatible with $H$, so that $\D{D}_*$ lifts to a
functor
$$\D{D}_*\colon \C{WCat}\To\C{qpm}.$$

The category
of square groups is a symmetric monoidal
category with the tensor product $\ca$ defined in \cite{qaI} that we now recall.

\begin{defn}\label{tens}
The \emph{tensor product} $M\ca N$ of square groups $M, N$ is
defined as follows. The group $(M\ca N)_e$ is generated by the
symbols $x\ca y$, $a\co b$ for $x\in M_e$, $y\in N_e$, $a\in M_{ee}$
and $b\in N_{ee}$, subject to the following relations
\begin{enumerate}
\item the symbol $a\co b$ is bilinear and central,
\item $x\ca (y_1+y_2)=x\ca y_1+ x\ca y_2+H(x)\co(y_2|y_1)_H$,
\item the symbol $x\ca y$ is left linear, $(x_1+x_2)\ca y=x_1\ca y+x_2\ca y$,
\item $P(a)\ca y=a\co (y|y)_H$.
\item $T(a)\co T(b)=-a\co b$,
\item $x\ca P(b)=\Delta(x)\co b$.
\end{enumerate}
The abelian group $(M\ca N)_{ee}$ is defined as the tensor product
$M_{ee}\otimes N_{ee}$. The homomorphism $$P\colon(M\ca N)_{ee}\To
(M\ca N)_{e}$$ is $P(a\otimes b)=a\co b$, and
$$H\colon(M\ca N)_{e}\To (M\ca N)_{ee}$$
is the unique quadratic map satisfying 
\begin{eqnarray*}
H(x\ca y)&=&\Delta(x)\otimes H(y)+ H(x)\otimes (y|y)_H,\\
H(a\co b)&=& a\otimes b-T(a)\otimes T(b),\\
(a\co b|-)_H&=&0,\\
(-|a\co b)_H&=&0,\\
(a\ca b|c\ca d)_H&=&(a|c)_H\otimes (b|d)_H.
\end{eqnarray*}
The unit for the tensor product is the square group $\Z_\ni$.
\end{defn}

\begin{thm}\label{prod}
Let $\C{C}\times\C{D}\r\C{E}\colon(A,B)\mapsto A\wedge B$ be a
biexact functor between Waldhausen categories. Then there are
morphisms of square groups
$$ \varphi^{ij}\colon
\D{D}^{sg}_i\C{C}\ca\D{D}^{sg}_j\C{D}\to\D{D}^{sg}_{i+j}\C{E},$$ for
$i$, $j$, $i+j\in\{0,1\}$, defined by
\begin{eqnarray*}
\varphi^{00}_e([A]\ca[C])&=&[A\wedge C],\\
\varphi^{01}_e([A]\ca[C\st{\sim}\r C'])&=&[A\wedge C\st{\sim}\r
A\wedge C'],\\
\varphi^{01}_e([A]\ca[C\rightarrowtail D\twoheadrightarrow
D/C])&=&[A\wedge C\rightarrowtail A\wedge D\twoheadrightarrow
A\wedge (D/C)],\\
\varphi^{10}_e([A\st{\sim}\r A']\ca[C])&=&[A\wedge C\st{\sim}\r
A'\wedge C],\\
\varphi^{10}_e([A\rightarrowtail B\twoheadrightarrow
B/A]\ca[C])&=&[A\wedge C\rightarrowtail B\wedge C\twoheadrightarrow
(B/A)\wedge C],\\
\varphi^{ij}_{ee}([A]\otimes[A']\otimes[C]\otimes[C'])&=&[A\wedge
C]\otimes[A'\wedge C'].
\end{eqnarray*}
such that the following diagram of square groups commutes
$$\xymatrix{&\D{D}^{sg}_1\C{C}\ca\D{D}^{sg}_1\C{D}\ar[rd]^{1\ca\partial}\ar[ld]_{\partial\ca1}&\\
\D{D}^{sg}_0\C{C}\ca\D{D}^{sg}_1\C{D}\ar[rd]^{\varphi^{01}}\ar@/_20pt/[rddd]_{1\ca\partial}&&
\D{D}^{sg}_1\C{C}\ca\D{D}^{sg}_0\C{D}\ar[ld]_{\varphi^{10}}\ar@/^20pt/[lddd]^{\partial\ca1}\\
&\D{D}^{sg}_1\C{E}\ar[d]^\partial&\\&\D{D}^{sg}_0\C{E}&\\&\D{D}^{sg}_0\C{C}\ca\D{D}^{sg}_0\C{D}\ar[u]_{\varphi^{00}}&}$$
\end{thm}

Now given a biexact functor
$\C{C}\times\C{D}\r\C{E}\colon(A,B)\mapsto A\wedge B$ we recover the
classical homomorphisms
\begin{eqnarray*}
\bar{\varphi}^{00}\colon K_0\C{C}\otimes
K_0\C{D}&\To&K_0\C{E},\\
\bar{\varphi}^{01}\colon
K_0\C{C}\otimes K_1\C{D}&\To&K_1\C{E},\\
\bar{\varphi}^{10}\colon K_1\C{C}\otimes K_0\C{D}&\To&K_1\C{E},
\end{eqnarray*}
from $\varphi^{ij}$ in Theorem \ref{prod} as follows. Given $i$,
$j$, $i+j\in\{0,1\}$,
\begin{eqnarray}\label{A}
\bar{\varphi}^{ij}(a\otimes b)&=&{\varphi}^{ij}_e(a\ca b).
\end{eqnarray}
Here we use the natural exact sequence
$$K_1\C{C}\hookrightarrow
\D{D}_1\C{C}\st{\partial}\To\D{D}_0\C{C}\twoheadrightarrow
K_0\C{C}$$ available for any Waldhausen category $\C{C}$ to identify
$K_1\C{C}$ with its image in $\D{D}_1\C{C}$, and we use the same
notation for an element in $\D{D}_0\C{C}$ and for its image in
$K_0\C{C}$. One can use the relations defining the tensor
product $\ca$ of square groups to check that the homomorphisms
$\bar{\varphi}^{ij}$ are well defined by the formula (\ref{A}) above.

In the proof of Theorem \ref{prod} we use a technical lemma about
square groups, which
measures the failure of $\ca$ to preserve certain coproducts.

Let
$M\cdot E$
be the $E$-fold
coproduct of a square group $M$, for any indexing set $E$.
We know from \cite{qaI} that
$\Z_\ni[E]=\Z_\ni\cdot E$, so we have canonical morphisms
$M\ca\iota_x\colon M\to M\ca \Z_\ni[E]$
for $x\in E$.
However, the natural comparison morphism
$$\iota=(M\ca\iota_x)_{x\in E}\colon M\cdot E\longrightarrow M\ca \Z_\ni[E]$$
is not an isomorphism.

Consider the homomorphisms
\begin{eqnarray*}
\overline{\Delta}\colon\Z[E]\otimes\Z[E]&\twoheadrightarrow&\Z[E],\\
H\colon\wedge^2\Z[E]&\r&
\ker\overline{\Delta}
,\\
q\colon
\ker\overline{\Delta}
&\twoheadrightarrow&\wedge^2\Z[E],
\end{eqnarray*}
where $\overline{\Delta}(e\otimes e)=e$ for $e\in E$ and
$\overline{\Delta}(e\otimes e')=0$ if $e\neq e'\in E$. Moreover,
$\wedge^2A$ is the exterior square of an abelian group $A$, i.e.
the quotient of $A\otimes A$ by the relations $a\otimes a=0$,
$a\in A$, $H(x\wedge y)=y\otimes x-x\otimes y$ and $q$ is induced
by the natural projection $A\otimes A\twoheadrightarrow\wedge^2
A\colon a\otimes b\mapsto a\wedge b$.
For any abelian group $A$ we consider 
$$(A\otimes\wedge^2\Z[E])^\otimes=\left(A\otimes\wedge^2\Z[E]
\mathop{\rightleftarrows}\limits^{A\otimes H}_{A\otimes q}
A\otimes\ker\overline{\Delta}
\right).$$
This is isomorphic to the square group defined in \cite[Section 1]{qaI}.
The construction is obviously functorial in $A$.

\begin{lem}\label{push}
For any square group $M$ and any set $E$ there is a pushout diagram in the category of square groups
$$\xymatrix{(\coker P\otimes \coker
P\otimes\wedge^2\Z[E])^\otimes\ar[r]\ar[d]\ar@{}[rd]|{\text{push}}
&M\cdot E\ar[d]^\iota\\
(M_{ee}\otimes\wedge^2\Z[E])^\otimes\ar[r]&M\ca\Z_\ni[E]\!\!\!\!\!\!\!\!\!\!\!\!\!}$$
which is natural in $M$ and $E$.
\end{lem}

\begin{proof}
One can check inductively by using \cite[Proposition 5]{qaI} and \cite[Section 5.6 (6)]{qaI} that there is
a map of central
extensions of square groups in the sense of \cite[Section 5.5]{qaI} as follows.
$$\xymatrix{(\coker P\otimes \coker P\otimes\wedge^2\Z[E])^\otimes\;
\ar@{^{(}->}[r]^-\mu
\ar[d]
&
M\cdot E\ar[d]_{\iota}\ar@{->>}[r]&\prod_EM\ar@{=}[d]\\
(M_{ee}\otimes\wedge^2\Z[E])^\otimes\;\ar@{^{(}->}[r]^-\nu&M\ca\Z_\ni[E]\ar@{->>}[r]&\prod_EM}$$
Here the left-hand morphism is induced by (\ref{qw}).
The morphism $\mu$ is completely determined by the homomorphism
$$\mu_{ee}\colon \coker P\otimes \coker
P\otimes\ker\overline{\Delta}\To(M\cdot E)_{ee}$$
defined by $\mu_{ee}(a\otimes b\otimes x\otimes
y)=P(\iota_x(a)|\iota_y(b))_H$ for $a,b\in \coker P$ and $x\neq y\in
E$.
Similarly $\nu$ is determined by the homomorphism
$$\nu_{ee}\colon M_{ee}\otimes\ker\overline{\Delta}\To M_{ee}\otimes \Z[E]\otimes\Z[E]$$
induced by the inclusion $\ker\overline{\Delta}\subset\Z[E]\otimes\Z[E]$.
It is straightforward to check that the square on the left is the desired pushout.
\end{proof}

\begin{proof}[Proof of Theorem \ref{prod}]
It is
obvious that
$\varphi^{ij}_{ee}$ is a well-defined abelian group homomorphism in all cases.
The square group morphism $\varphi^{00}$ is well-defined as a consequence of \cite[Proposition 34]{qaI}.

Let $E$ be the set of objects of $\C{D}$, excluding
$*$, so that
$\Z_\ni[E]=\D{D}_0^{sg}\C{D}$, and let $M=\D{D}_1^{sg}\C{C}$.
To see that $\varphi^{10}$ is well defined by the formulas
in the statement we note
that it is just the morphism determined, using
Lemma \ref{push}, by the square group morphisms
$$
(M_{ee}\otimes\wedge^2\Z[E])^\otimes
\st\xi\longrightarrow
\D{D}^{sg}_1\C{E}
\st\zeta\longleftarrow
M\cdot E
$$
defined as follows.
The square group morphism $\xi$ is completely determined by
$$\xi_{ee}=\varphi^{10}_{ee}\colon(\D{D}_0\C{C})^\abb\otimes(\D{D}_0\C{C})^\abb
\otimes\ker\overline{\Delta}\To(\D{D}_0\C{E})^\abb\otimes(\D{D}_0\C{E})^\abb.$$
For each $A\in E$,
the component $\zeta\circ \iota_A:\D{D}_1^{sg}\C{C}\to \D{D}_1^{sg}\C{C}\cdot E\to \D{D}^{sg}_1\C{E}$
is the unique square group morphism such that
\begin{equation*}
\xymatrix@C=60pt{
\D{D}^{sg}_1\C{C}\ar[r]^{\iota_A}
\ar[d]_\partial&\D{D}^{sg}_1\C{C}\ar[r]^{\zeta}\cdot E
&\D{D}^{sg}_1\C{E}\ar[d]^\partial\\
\D{D}^{sg}_0\C{C}\ar[r]^{\iota_A}
&\D{D}^{sg}_0\C{C}\ar[r]^{\varphi^{00}}\cdot E
&\D{D}^{sg}_1\C{E}
}\end{equation*}
coincides with
the morphism of quadratic pair modules
$$\D{D}_*(\,\cdot\,\wedge A)\colon\D{D}_*\C{C}\r\D{D}_*\C{E}$$
induced by the exact functor $\,\cdot\,\wedge A\colon\C{D}\r\C{E}$.

By using this alternative definition of $\varphi^{10}$ in terms of Lemma \ref{push} it is also immediate that the lower right cell in the diagram of the statement is commutative.

To see that $\varphi^{01}$ is well-defined and that the lower left cell of the diagram in the
statement commutes one proceeds similarly, using the fact that $\ca$ is symmetric.

Now we just need to check that the upper cell is commutative. For
this it is enough to show that the following equalities hold

\begin{scriptsize}
\begin{eqnarray*}
\varphi^{01}((\partial[A\st{\sim}\r A'])\ca[C\st{\sim}\r
C'])&=&\varphi^{10}([A\st{\sim}\r A']\ca(\partial[C\st{\sim}\r
C'])),\\
\varphi^{01}((\partial[A\st{\sim}\r A'])\ca[C\rightarrowtail
D\twoheadrightarrow D/C])&=&\varphi^{10}([A\st{\sim}\r
A']\ca(\partial[C\rightarrowtail D\twoheadrightarrow D/C])),\\
\varphi^{01}((\partial[A\rightarrowtail B\twoheadrightarrow
B/A])\ca[C\st{\sim}\r C'])&=&\varphi^{10}([A\rightarrowtail
B\twoheadrightarrow B/A]\ca(\partial[C\st{\sim}\r
C'])),\\
\varphi^{01}((\partial[A\rightarrowtail B\twoheadrightarrow
B/A])\ca[C\rightarrowtail D\twoheadrightarrow
D/C])&=&\varphi^{10}([A\rightarrowtail B\twoheadrightarrow
B/A]\ca(\partial[C\rightarrowtail D\twoheadrightarrow D/C])).
\end{eqnarray*}%
\end{scriptsize}%
This is a tedious but straightforward task which makes use of the
laws of stable quadratic modules and the tensor product of square
groups, the elementary properties of a biexact functor, and the relations (1),
(2), (6), (7) and (8) in Definition \ref{LA}.
\end{proof}

\section{Natural transformations and induced homotopies on $\D{D}_*$}

\noindent
In this section we define induced homotopies along the functor $\D{D}_*$ from section \ref{1}.

\begin{defn}\label{homoto}
Two morphisms $f,g\colon C_*\r D_*$ are \emph{homotopic} $f\simeq g$
if there exists a function $\alpha\colon C_0\r D_1$ satisfying
\begin{itemize}
\item
$\alpha(c_0+d_0)=\alpha(c_0)+\alpha(d_0)+\grupo{f_0(d_0),-f_0(c_0)+g_0(c_0)}_{D_*}$,
\item $g_0(c_0)=f_0(c_0)+\partial_{D_*}\alpha(c_0)$,
\item $g_1(c_1)=f_1(c_1)+\alpha\partial_{C_*}(c_1)$.
\end{itemize}
Such a function is called a \emph{homotopy} $\alpha\colon f\rr g$,
also denoted by
\begin{equation*}
\xymatrix@C=35pt{C_*\ar@/^10pt/[r]^f_{\;}="a"\ar@/_10pt/[r]_{g}^{\;}="b"&
D_*.\ar@{=>}"a";"b"^{\alpha}}
\end{equation*}
\end{defn}

The category $\C{squad}$ of stable quadratic modules is a
$2$-category with $2$-morphisms given by homotopies. It is indeed
a category enriched over groupoids, see \cite[Proposition
7.2]{2hg1}. The vertical composition of $2$-morphisms
$$\xymatrix@C=40pt{C_*\ar@/^20pt/[r]^f_{\;}="a"\ar[r]|{g}^{\;}="b"_{\;}="c"\ar@/_20pt/[r]_h^{\;}="d"&D_*\ar@{=>}"a";"b"^\alpha\ar@{=>}"c";"d"^\beta}$$
denoted by $\beta\vc\alpha\colon f\rr h$, is defined as
$(\beta\vc\alpha)(x)=\alpha(x)+\beta(x)$ for $x\in C_0$. The
identity $2$-morphism
\begin{equation*}
\xymatrix@C=35pt{C_*\ar@/^10pt/[r]^f_{\;}="a"\ar@/_10pt/[r]_{f}^{\;}="b"&
C_*.\ar@{=>}"a";"b"^{1_f}}
\end{equation*}
is clearly given by the function $C_0\r C_1\colon x\mapsto 0$. The
horizontal compositions $h\alpha\colon hf\rr hg$ and $\alpha k\colon fk\rr gk$
in the diagram
\begin{equation*}
\xymatrix@C=35pt{C_*\ar[r]^{k}&D_*\ar@/_15pt/[r]_{g}^{\;}="a"
\ar@/^15pt/[r]^{f}_{\;}="b"&E_*
\ar[r]^{h}&K_*,\ar@{=>}"b";"a"^\alpha}
\end{equation*}
are
defined by the functions $h_1\alpha\colon D_0\r K_1$ and $\alpha
k_0\colon C_0\r E_1$, respectively.

The homotopy category obtained by imposing the homotopy relation on the
full subcategory $\C{squad}_f\subset\C{squad}$
given by $0$-free objects is equivalent
to the homotopy category of stable quadratic modules
$$\C{squad}_f/\!\simeq\;\st{\sim}\To\;\ho\C{squad},$$
compare \cite[Proposition 7.7]{2hg1}.

The category $\C{WCat}$ of Waldhausen categories and exact functors is also a $2$-category, where
a $2$-morphism $\eps\colon F\rr G$ between two exact functors $F,G\colon\C{C}\r\C{D}$ is a natural
transformation $\eps$ given by weak equivalences $\eps(A)\colon F(A)\st{\sim}\r G(A)$ in $\C{D}$ for any object
$A$ in $\C{C}$.

\begin{thm}\label{2}
The construction $\D{D}_*\colon\C{WCat}\r\C{squad}$ defines a $2$-functor.
\end{thm}

\begin{proof}
The homotopy $\D{D}_*\eps\colon\D{D}_0\C{C}\r\D{D}_1\C{D}$ induced
by a $2$-morphism $\eps\colon F\rr G$ in $\C{WCat}$ is defined on
generators $[A]\in\D{D}_0\C{C}$ by the formula
$$(\D{D}_*\eps)([A])\;=\;{}-[\eps(A)\colon F(A)\st{\sim}\r G(A)],$$
and then extended to the whole group $\D{D}_0\C{C}$ by using the
first equation in Definition \ref{homoto}. This is a well-defined
homotopy, compare \cite[IV.4.5]{ch4c}. Now we have to check that
$\D{D}_*$, defined in this way, preserves the vertical composition
of $2$-morphisms, the identity $2$-morphisms, and the horizontal
composition of a $1$-morphism and a $2$-morphism.

Let
$$\xymatrix@C=40pt{\C{C}\ar@/^20pt/[r]^F_{\;}="a"\ar[r]|{G}^{\;}="b"_{\;}="c"\ar@/_20pt/[r]_H^{\;}="d"&\C{D}\ar@{=>}"a";"b"^\varepsilon\ar@{=>}"c";"d"^\delta}$$
be a diagram of vertically composable $2$-morphisms in $\C{WCat}$.
Given an object $A$ in $\C{C}$ we have
\begin{eqnarray*}
((\D{D}_*\delta)\vc(\D{D}_*\varepsilon))([A])&=&(\D{D}_*\varepsilon)([A])+(\D{D}_*\delta)([A])\\
&=&-[\eps(A)\colon F(A)\st{\sim}\r G(A)]-[\delta(A)\colon
G(A)\st{\sim}\r H(A)]\\
\text{{\small by Definition \ref{LA}
(6)}}\qquad&=&-[(\delta\vc\varepsilon)(A)=\delta(A)\eps(A)\colon
F(A)\st{\sim}\r G(A)]\\
&=&(\D{D}_*(\delta\vc\varepsilon))([A]),
\end{eqnarray*}
so vertical composition is preserved.

By Definition \ref{LA} (4) $\D{D}_*$ takes the identity natural
transformation to the identity homotopy.

Finally given a diagram in $\C{WCat}$ and objects $A$ in $\C{D}$
and $B$ in $\C{C}$
\begin{equation*}
\xymatrix@C=35pt{\C{C}\ar[r]^{K}&\C{D}\ar@/_15pt/[r]_{G}^{\;}="a"
\ar@/^15pt/[r]^{F}_{\;}="b"&\C{E}
\ar[r]^{H}&\C{K}\ar@{=>}"b";"a"^\varepsilon}
\end{equation*}
we have
\begin{eqnarray*}
(\D{D}_1H)(\D{D}_*\varepsilon)([A])&=&(\D{D}_1H)(-[\eps(A)\colon
F(A)\st{\sim}\r G(A)])\\
&=&-[(H\eps)(A)=H(\eps(A))\colon HF(A)\st{\sim}\r HG(A)]\\
&=&(\D{D}_*(H\varepsilon))([A]),\\
(\D{D}_*\varepsilon)(\D{D}_0K)([B])&=&(\D{D}_*\varepsilon)([K(B)])\\
&=&-[(\eps K)(B)=\eps (K(B))\colon FK(B)\st{\sim}\r GK(B)]\\
&=&(\D{D}_*(\varepsilon K))([B]).
\end{eqnarray*}
Now the proof is finished.
\end{proof}

\begin{rem}
The homotopies defined in Theorem \ref{2} are constructed by using just one kind of degree $1$ generators of
$\D{D}_*$, namely those given by weak equivalences. In case we have a cofiber sequence $F\rightarrowtail
G\twoheadrightarrow H$ of exact functors $F,G,H\colon \C{C}\r\C{D}$ one can define a homotopy
using the other class of degree $1$ generators, given by cofiber sequences, to
give a direct proof of the additivity theorem \cite[Proposition 1.3.2 (4)]{akts} for the algebraic
model of the $1$-type $\D{D}_*$.
\end{rem}

\section{Proof of Theorem \ref{main}}

In this section we use the notions of crossed module and crossed
complex in the category of groups and in the category of
groupoids. There are different but equivalent ways of presenting
these objects 
depending on a
series of conventions such as using left or right actions, choice
of basepoint of an $n$-simplex, etc. In this paper we adopt the
conventions which are compatible with \cite{EZcc}. As examples of
a crossed complex we can mention the fundamental crossed complex
$\pi_{CW}Y$ of a $CW$-complex $Y$ and the fundamental crossed
complex $\pi X$ of a simplicial set $X$; these are related by the
natural identification $\pi X=\pi_{CW}\abs{X}$ where
$\abs{\,\cdot\,}$ denotes the geometric realization functor from
simplicial sets to $CW$-complexes. See \cite{EZcc} for further
details and references. The reader who is unfamiliar with
simplicial techniques is referred to the texts~\cite{sht,saat}.

\begin{defn}\label{cmes}
A \emph{crossed module (of groups)} is a group homomorphism $\partial\colon M\r N$ such that $N$ acts (on the right)
on $M$ and the following equations are satisfied for $m,m'\in M$ and $n\in N$.
\begin{eqnarray}
\partial(m^n)&=&-n+\partial(m)+n,\label{cm1}\\
m^{\partial(m')}&=&-m'+m+m'.\label{cm2}
\end{eqnarray}

Morphisms of crossed modules are defined by commutative squares of group homomorphisms which are compatible with
the actions in the obvious way. Such a morphism is a \emph{weak
  equivalence} if it induces isomorphisms between the
kernels and cokernels of the homomorphisms $\partial$.
\end{defn}

A {\em crossed complex (of groups)} $(C,\partial)$ is given by groups and homomorphisms
$$\cdots\r C_n\st{\partial_n}\To C_{n-1}\r\cdots\r C_3 \st{\partial_3}\To C_2\st{\partial_2}\To C_1,\qquad \partial^2=0,$$
where $\partial_2\colon C_2\to C_1$ is a crossed module as above,
$\{C_n,\; n\ge3\}$ is a chain complex of modules over
$\coker[\partial_2\colon C_2\to C_1]$, and the homomorphism
$\partial_3\colon C_3\to C_2$ is compatible with the actions of
$C_1$. Note that the kernel of $\partial_2$ (and in particular the
image of $\partial_3$) is central in $C_2$ by \eqref{cm2}.

We sometimes adopt the convention $C_0=\{*\}$
for crossed complexes of groups, since they are the `one object'
cases of the more general crossed complexes of groupoids.
In the general case the groupoid $C_1$ has source and target functions
$s,t:C_1\to C_0$,
while for $n\geq 2$ the groupoids $C_n$ are {\em totally disconnected},
that is, $s=t\colon C_n\to C_0$.
The action of $C_1$ is now only partially defined, given by functions $C_n(y,y)\times C_1(x,y)\to C_n(x,x)$ for objects $x,y\in C_0$.

If $Y$ is a $CW$-complex the crossed complex of groupoids
$\pi_{CW}Y$ is defined by the relative homotopy groupoids of the
skeletons $Y^n$ of $Y$ based at the set of points $Y^0$
\begin{equation*}
\begin{array}{r}\cdots\r \pi_n(Y^n,Y^{n-1},Y^0)\st{\partial_n}\To
\pi_{n-1}(Y^{n-1},Y^{n-2},Y^0)\r\cdots\\\cdots\r
\pi_2(Y^2,Y^1,Y^0)\st{\partial_2}\To \pi_1(Y^1,Y^0).
\end{array}
\end{equation*}
This is a crossed complex of groups if and only if $Y$ is a
reduced $CW$-complex, i.e. $Y^0=\set{*}$.

If $X$ is a simplicial set, then the crossed complex of groupoids
$\pi X$ has object set $X_0$ and an explicit presentation by
generators $x\in(\pi X)_k$ for each $k$-simplex $x$ of $X$, with
source ${d_1}^{\!\!k}x$ (and target $d_0x$ if $k=1$), and
relations given by the boundaries
\begin{eqnarray}
\label{dsimp}
\partial_2x&=&{}-d_1x+d_0x+d_2x\\
\nonumber
\partial_3x&=&d_2x+d_0x^{\,{d_2}^{\!\!2}x}-d_3x-d_1x\\
\nonumber
\partial_kx&=&d_0x^{\,{d_2}^{\!\!k-1}x}+\sum_{i=1}^r(-1)^id_ix
\end{eqnarray}
together with the relation $x=0$ in $\pi X$ whenever
$x$ is a degenerate simplex of $X$;
compare \cite[Example 1.2]{EZcc}.

The category $\C{ccplx}$ of crossed complexes is symmetric monoidal
with respect to a tensor product $\otimes$, introduced in \cite{tpcc}.
Suppose $C$ and $D$ are crossed complexes of groupoids. Then, using the
conventions of \cite[Definition 1.4]{EZcc}, the tensor product $C\otimes D$
has a presentation with generators $c\otimes d\in(C\otimes D)_{m+n}$ for each $c\in C_m$ and $d\in D_n$, where $c\otimes d$ has source $sc\otimes sd\in (C\otimes D)_0$ (and target $tc\otimes d$ or $c\otimes td$ if $(m,n)=(1,0)$ or $(0,1)$).
The relations are
\begin{align}
 \label{bilin3}
(c\otimes d)^{a\otimes b}&=\left\{
\begin{array}{cl}
\qquad c^a\otimes d^{\phantom b} \qquad
                          & \textrm{ if }b=sd \textrm{ and } m\geq 2,\\
\qquad c^{\phantom a}\otimes d^b \qquad
                          & \textrm{ if }a=sc \textrm{ and } n\geq 2,
\end{array}\right.
\\
\nonumber
c\otimes (d+d')&=\left\{
\begin{array}{cl}
\quad c\otimes d\quad \;\; + \;\; c\otimes d'\quad
                          & \textrm{ if }n\geq2 \textrm{ or } m=0,\\
(c\otimes d)^{sc\otimes d'} + \;\; c\otimes d'
                          & \textrm{ if }n=1 \textrm{ and } m\geq 1,
\end{array}\right.
\\
\nonumber
(c+c')\otimes d&=\left\{
\begin{array}{cl}
c\otimes d \;\; + \;\; c'\otimes d
                          & \textrm{ if }m\geq2 \textrm{ or } n=0,\\
\qquad c'\otimes d\;\;+\;\;(c\otimes d)^{c'\otimes sd}
                          & \textrm{ if }m=1 \textrm{ and } n\geq 1,
\end{array}\right.
\\
\nonumber
\partial_{m+n}
(c\otimes d)&=\left\{
\begin{array}{cl}
-sc\otimes d-c\otimes td+tc\otimes d+c\otimes sd
                          & \textrm{ if }m,n=1,\\
\partial_m\I (c\otimes d)\;\;+\;\;(-1)^m\;\partial_n\II(c\otimes d)
                          & \textrm{ otherwise.}
\end{array}\right.
\end{align}
The symbols $\partial_k\I(c\otimes d)$ and $\partial_k\II(c\otimes d)$ are defined
for $k\geq2$ by $\partial_k c\otimes d$ and   $c\otimes \partial_k d$ respectively,
for $k=1$ by $-{sc\otimes d}+(tc\otimes d)^{c\otimes sd}$
and $-{c\otimes sd}+(c\otimes td)^{sc\otimes d}$,
and vanish for $k=0$.

This tensor product satisfies the following crucial property:
given two $CW$-complexes $Y$, $Z$ there is a natural isomorphism
$\pi_{CW}Y\otimes\pi_{CW}Z\cong\pi_{CW}(Y\times Z)$ \cite[Theorem
3.1 (iv)]{cscc} satisfying the usual coherence properties. For
products of simplicial sets there is an Eilenberg-Zilber theorem
proved in \cite{EZcc}. As examples of monoids in the category of
crossed complexes we can cite the fundamental crossed complex
$\pi_{CW}M$ of a $CW$-monoid $M$, and the crossed cobar
construction $\underline{\Omega}_{\C{Crs}}X$ on a $1$-reduced
simplicial set $X$, see \cite{tcc}. As a consequence of
\cite{EZcc} the fundamental crossed complex $\pi N$ of a
simplicial monoid $N$ is also a monoid in $\C{ccplx}$.

For our purposes it will be convenient to have a small model for the fundamental crossed complex of the diagonal
of a bisimplicial set. This is achieved by the following definition.

\begin{defn}
The {\em total crossed complex} $\Pi(X)$ of a bisimplicial set $X$ is
the coend
\begin{eqnarray}\label{bigpi}
\Pi(X)&=&\int^{m,n}\pi(\Delta[m])\otimes\pi(\Delta[n])\;\cdot\;X_{m,n}
\;.\end{eqnarray}
Here $\Delta[k]$ is the $k$-simplex, $k\geq 0$, and $C\cdot E$ is the $E$-fold coproduct of a crossed complex
$C$ over an indexing
set $E$; see \cite[IX.6]{cwm} for more details on coend calculus.
\end{defn}

The following lemma gives an explicit presentation
in terms of
generators and relations.

\begin{lem}\label{presen}
Suppose $X$ is a horizontally-reduced bisimplicial set,
in the sense that $X_{0,*}=\Delta[0]$.
Then $\Pi(X)$ is the crossed complex of groups with one generator $x_{m,n}$ in $\Pi(X)_{m+n}$ for each
$x_{m,n}\in X_{m,n}$ and subject to the following
relations:
\begin{eqnarray*}
x_{m,n}&=&0\quad\text{if $x_{m,n}$ is degenerate in $X_{m,n}$}
\;,\\
\partial_2x_{1,1}
&=&{}
-d_0^vx_{1,1}
+d_1^vx_{1,1}
\;,\\
\partial_2x_{2,0}
&=&{}
-d_1^hx_{2,0}+d_0^hx_{2,0}+d_2^hx_{2,0}
\;,\\
\partial_3x_{1,2}
&=&{}
-d_2^vx_{1,2}-d_0^vx_{1,2}+d_1^vx_{1,2}
\;,\\
\partial_3x_{2,1}
&=&{}
 d_2^hx_{2,1}
+{d_0^hx_{2,1}}^{d_2^hd_1^vx_{2,1}}
-d_1^vx_{2,1}
-d_1^hx_{2,1}
+d_0^vx_{2,1}
\;,\\
\partial_3x_{3,0}
&=&{}
d^h_2x_{3,0} + {d^h_0x_{3,0}}^{(d_2^h)^2x_{3,0}} -d^h_3x_{3,0} -d^h_1x_{3,0}
\;.\end{eqnarray*}
For $m\ge1$ and $m+n\ge4$, the boundary relations are abelian:
\begin{eqnarray*}
\partial_{m+n}x_{m,n}
&=&{}
{d^h_0x_{m,n}}^{(d_2^h)^{m-1}(d_1^v)^nx_{m,n}} +\sum_{i=1}^m(-1)^i d^h_ix_{m,n}
+\sum_{j=0}^n(-1)^{m+j} d^v_jx_{m,n}.
\end{eqnarray*}
The last summation is trivial if $n=0$;
all the other terms are trivial if $m=1$.
\end{lem}

\begin{proof}
The proof follows from the abstract coend definition by using the presentations for $\pi\Delta[k]$ and for the tensor product of crossed complexes given above, and simplifying via some straightforward calculations in crossed complexes. We give the full details.

The total crossed complex $\Pi(X)$ has generators in dimension
$m+n$ given by
$$(c_m\otimes d_n;\;x_{m',n'}),\qquad
\text{ for }\,c_m\in(\pi\Delta[m'])_m,\;d_n\in(\pi\Delta[n'])_n,\;x_{m',n'}\in X_{m',n'},$$
subject to the relations listed in 
\eqref{bilin3} and the coend identifications
$$(\sigma_*c_m\otimes \tau_*d_n;\,x_{m'',n''})
=(c_m\otimes d_n;\,(\sigma,\tau)^*x_{m'',n''}),
\text{ for }
(\sigma,\tau)\in\Delta(m',m'')\times\Delta(n',n'').$$
By the first three relations of \eqref{bilin3} we need only consider those generators for which $c_m$ is a generator of $\pi\Delta[m']$, given by an $m$-simplex of $\Delta[m']$, and similarly for $d$. Since $\Delta[m']_m=\Delta(m,m')$ we have the coend identification
\begin{align*}
(c_m\otimes d_n;\;x_{m',n'})&=(\delta_m\otimes\delta_n;\;(c_m,d_n)^*x_{m',n'})
\end{align*}
and we can assume $m=m'$, $n=n'$ and $c$, $d$ are top-dimensional simplices $\delta$ given by the `identity' simplicial operators $1$.

Generators of $\Pi(X)$ are thus identified simply with elements $x_{m,n}\in X_{m,n}$, and $\Pi(X)$ is a crossed complex of groups, $(\Pi(X))_0=X_{0,0}=\{*\}$.
The horizontally or vertically degenerate elements $x_{m,n}$ are zero,
since the degeneracy may be transferred to the generator
of $\pi\Delta[m]$ or $\pi\Delta[n]$ respectively.
In particular
we have $x_{0,n}=0$  whenever $n>0$, since $X$ is horizontally reduced.

The boundary maps on $\pi\Delta[m]$ and $\pi\Delta[n]$ are specified by the
relations \eqref{dsimp}.
Transferring them to horizontal and vertical face maps on $x_{m,n}$, the
fourth relation of \eqref{bilin3} gives the following when $m+n\le3$
\begin{align*}
(\partial_{1+1}&(\delta_1\otimes\delta_1);
\;x_{1,1})=
(-d_1\delta_1\otimes\delta_1-\delta_1\otimes d_0\delta_1
 +d_0\delta_1\otimes\delta_1+\delta_1\otimes d_1\delta_1;\;x_{1,1})\\
&={}-0\;-\;(\delta_1\otimes\delta_1;\;d_0^vx_{1,1})
\;+\;0\;+\;(\delta_1\otimes\delta_1;\;d_1^vx_{1,1})
\,,\\
(\partial_{2+0}&(\delta_2\otimes\delta_0);
\;x_{2,0})=((-d_1\delta_2+d_0\delta_2+d_2\delta_2)\otimes\delta_0;\;x_{2,0})\\
&=-(\delta_2\otimes\delta_0;\;d_1^hx_{2,0})
  +(\delta_2\otimes\delta_0;\;d_0^hx_{2,0})
  +(\delta_2\otimes\delta_0;\;d_2^hx_{2,0})
\,,
\end{align*}
\begin{align*}
(\partial_{1+2}&(\delta_1\otimes\delta_2);
\;x_{1,2})=
(\partial_1\I (\delta_1\otimes\delta_2)
-\partial_2\II(\delta_1\otimes\delta_2);\; x_{1,2})\\
&=(-d_1\delta_1\otimes\delta_2
 +{(d_0\delta_1\otimes\delta_2)}
     ^{\delta_1\otimes d_1\delta_2}
      -\delta_1\otimes(-d_1\delta_2+d_0\delta_2+d_2\delta_2);\;x_{1,2})\\
&={}
-0+0 -(\delta_1\otimes\delta_2;\;d_2^v x_{1,2})
     -(\delta_1\otimes\delta_2;\;d_0^v x_{1,2})
     +(\delta_1\otimes\delta_2;\;d_1^v x_{1,2})
\,,\\
(\partial_{2+1}&(\delta_2\otimes\delta_1)
;
\;x_{2,1})=
(\partial_2\I(\delta_2\otimes\delta_1)+\partial_1\II(\delta_2\otimes\delta_1)
;\;x_{2,1})
\\&=
((-d_1\delta_2+d_0\delta_2+d_2\delta_2)\otimes\delta_1-\delta_2\otimes d_1\delta_1
+\delta_2\otimes d_0\delta_1;\;x_{2,1})
\\&=
(d_2\delta_2\otimes\delta_1+(d_0\delta_2\otimes\delta_1)^{d_2\delta_2\otimes d_1\delta_1}
\\&
\qquad
-(d_1\delta_2\otimes\delta_1)^{\partial_2\delta_2\otimes d_1\delta_1}-\delta_2\otimes d_1\delta_1+\delta_2\otimes d_0\delta_1;\;x_{2,1})
\\&=
 (\delta_2\otimes\delta_1;\;   d_2^hx_{2,1})
+(\delta_2\otimes\delta_1;\;   d_0^hx_{2,1})^{(\delta_2\otimes\delta_1;\;d_2^hd_1^vx_{2,1})}
\\&
\qquad
-(\delta_2\otimes\delta_1;\;  d_1^vx_{2,1})
-(\delta_2\otimes\delta_1;\;  d_1h^x_{2,1})
+(\delta_2\otimes\delta_1;\;  d_0^vx_{2,1})
\,,\\
(\partial_{3+0}&(\delta_3\otimes\delta_0)
;
\;x_{3,0})=
((d_2\delta_3+(d_0\delta_3)^{{d_2}^2\delta_3}-d_3\delta_3-d_1\delta_3)\otimes\delta_0
;\;x_{2,0})\\
&=
(\delta_3\otimes\delta_0;\;d_2^hx_{3,0})+
(\delta_3\otimes\delta_0;\;d_0^hx_{3,0})^{(\delta_3\otimes\delta_0;\,{d_2^h}^2x_{3,0})}
\\&\qquad-(\delta_3\otimes\delta_0;\;d_3^hx_{3,0})-(\delta_3\otimes\delta_0;\;d_1^hx_{3,0})
.\end{align*}
These correspond exactly to the first five boundary relations claimed in the statement
of this Lemma.
Note that in the penultimate relation we used the fact that $X$ is horizontally reduced, so the action disappears in the expansion of $\partial_1\II$, and we also
used the third relation of \eqref{bilin3} and the crossed module axiom \eqref{cm2} in the expansion of $\partial_2\I$.

The remaining relation to prove, in the chain complex $(\Pi(X))_{\geq3}$, is a combination of the following cases:
\begin{align*}
(\partial_{m+0}&(\delta_m\otimes\delta_0);\;x_{m,0})
=
\biggl(\,(d_0\delta_m\otimes\delta_0)^{{d_2}^{m-1}\delta_m\otimes\delta_0}+
\sum_{i=1}^m(-1)^i\,d_i\delta_m\otimes\delta_0;\;x_{m,0}\biggr)
\\
&
=
(\delta_m\otimes\delta_0;\;d_0^hx_{m,0})^{(\delta_m\otimes\delta_0;\,(d_2^h)^{m-1}x_{m,0})}+\sum_{i=1}^m(-1)^i
(\delta_m\otimes\delta_0;\;d_i^hx_{m,0})
\,,
\\
(\partial_{m+1}&(\delta_m\otimes\delta_1);\;x_{m,1})
=(\partial_m\I(\delta_m\otimes\delta_1)+(-1)^m\partial_1\II(\delta_m\otimes\delta_1);\;x_{m,1})
\\&=
\biggl((d_0\delta_m\otimes\delta_1)^{{d_2}^{m-1}\delta_m\otimes d_1\delta_1}+\sum_{i=1}^m(-1)^i\,d_i\delta_m\otimes\delta_1
\\&\qquad\qquad\qquad
+(-1)^m(-\delta_m\otimes d_1\delta_1+\delta_m\otimes d_0\delta_1
);\;x_{m,1}\biggr)
\\&=
(\delta_m\otimes\delta_1;\;d_0^hx_{m,1})^{(\delta_m\otimes\delta_1
;\,{d_2^h}^{m-1}d_1^vx_{m,1})}
+\sum_{i=1}^m(-1)^i\,(\delta_m\otimes\delta_1;\;d_i^hx_{m,1})
\\&\qquad\qquad\qquad
+(-1)^m\bigl(
-(\delta_m\otimes\delta_1;\;d_1^vx_{m,1})
+(\delta_m\otimes\delta_1;\;d_0^vx_{m,1})
\bigr)
\,,
\\
(\partial_{1+n}&(\delta_1\otimes\delta_n);\;x_{1,n})
=\bigl(0-\partial_n\II(\delta_1\otimes\delta_n);\;x_{1,n}\bigr)
=-\biggl(\sum_{j=0}^n(-1)^j\delta_1\otimes d_j\delta_n;\;x_{1,n}\biggr)
\\
&
=\sum_{j=0}^n(-1)^{1+j}(\delta_1\otimes\delta_n;\;d_j^vx_{1,n})
\,,
\end{align*}
\begin{align*}
(\partial_{m+n}&(\delta_m\otimes\delta_n);\;x_{m,n})
=(\partial_m\I(\delta_m\otimes\delta_n)+(-1)^m\partial_1\II(\delta_m\otimes\delta_n);\;x_{m,n})
\\
&
=\biggl((d_0\delta_m\otimes\delta_n)^{{d_2}^{m-1}\delta_m\otimes d_1^n\delta_n}+\sum_{i=1}^m(-1)^i\,d_i\delta_m\otimes\delta_n
\\&\qquad\qquad\qquad\qquad
\;+\;
(-1)^m\;\sum_{j=0}^n(-1)^{j}\delta_m\otimes d_j\delta_n;\;\;x_{m,n}\biggr)
\\
&
=(\delta_m\otimes\delta_n;\;d_0^nx_{m,n})^{
(\delta_m\otimes\delta_n;\,{d_2^h}^{m-1}{d_1^v}^nx_{m,n})}+
\sum_{i=1}^m(-1)^i(\delta_m\otimes\delta_n;\;d_i^hx_{m,n})
\\&\qquad\qquad\qquad
+
\sum_{j=0}^n(-1)^{m+j}(\delta_m\otimes\delta_n;\;d_j^vx_{m,n})
.
\end{align*}
\end{proof}

The following results are natural generalizations of the Eilenberg--Zilber theorem for
crossed complexes given in~\cite{EZcc}.

\begin{thm}\label{biEZ}
There is a natural homotopy equivalence (in fact, a strong deformation
retraction) of crossed complexes between the total crossed complex of
a bisimplicial set $X$ and the fundamental crossed complex of its
diagonal,
\[\EZDIAG{\pi\diag(X)}{\Pi(X).}{a'}{b'}{\phi'}\]
\end{thm}

\begin{proof}
As observed for example in~\cite[Proposition B.1]{htgamma},
the diagonal of a bisimplicial set $X$ may be expressed as a coend
\[
\diag(X)\;\;\cong\;\;\int^m\Delta[m]\times X_{m,*}.
\]
Since each $X_{m,*}$ is the coend of $\Delta[n]\cdot X_{m,n}$,
and $\pi$ preserves colimits,
\begin{eqnarray*}
\pi\diag(X)&\cong&\pi\int^{m,n}\Delta[m]\times\Delta[n]\;\cdot\; X_{m,n}
\\&\cong&\int^{m,n}\pi(\Delta[m]\times\Delta[n])
\;\cdot\; X_{m,n}
\;.
\end{eqnarray*}
The result therefore follows from the Eilenberg--Zilber equivalence
\[\EZDIAG{\pi(\Delta[m]\times\Delta[n])}{\pi\Delta[m]\otimes\pi\Delta[n]}{a}{b}{\phi}\]
given in
\cite[Theorem 3.1]{EZcc} (see also
\cite[Section 3]{smcsecc}).
\end{proof}

\begin{thm}\label{laxmon}
Given two
bisimplicial sets $X$, $Y$, there is a natural deformation retraction
\begin{equation}
\EZDIAG{\Pi(X\times Y)}{\Pi(X)\otimes \Pi(Y).}{a''}{b''}{\phi''}
\end{equation}
Moreover, the following
diagram of `shuffle maps' commutes:
\begin{equation}\label{bdiag}{\xymatrix{
\Pi X\otimes \Pi Y\ar[rrr]^(0.4){b'\otimes b'} \ar[d]_{b''}&&&
\pi\diag X\otimes \pi\diag Y\ar[d]^b\\
\Pi(X\times Y)\ar[rr]^-{b'}&&
\pi\diag(X\times Y)\ar[r]^-\cong&
\pi(\diag X\times\diag Y)
}}\end{equation}
\end{thm}

\begin{proof}
The natural homotopy equivalence of the objects
\begin{eqnarray*}
\Pi(X\times Y)&\cong&\int^{p,p',q,q'}\pi(\Delta[p]\times\Delta[p'])
\otimes \pi(\Delta[q]\times\Delta[q'])\cdot X_{p,q}\times Y_{p',q'},\\
\Pi(X)\otimes\Pi(Y)&\cong&\int^{p,p',q,q'}\pi\Delta[p]\otimes
\pi\Delta[q]\otimes
\pi\Delta[p']\otimes
\pi\Delta[q']\cdot X_{p,q}\times Y_{p',q'},\\
\end{eqnarray*}
is defined using the symmetry $\pi\Delta[q]\otimes\pi\Delta[p']\cong\pi\Delta[p']\otimes\pi\Delta[q]$
and the Eilenberg--Zilber equivalence, see \cite{EZcc}. The commutativity of the
diagram (\ref{bdiag}) follows from standard properties of the shuffle map.
\end{proof}

\begin{exm}\label{ll}
Suppose $X$, $Y$ are bisimplicial sets, with $x\in X_{1,0}$ and $y\in
Y_{1,0}$ corresponding to generators in degree one of $\Pi X$ and $\Pi
Y$ respectively.
Then by \cite[2.6]{EZcc} we have $b''(x\otimes y)
\in \Pi(X\times Y)_2$ given by
\begin{eqnarray*}
b''(x\otimes y)&=&-(s_0^hx,s_1^hy)+(s_1^hx,s_0^hy).
\end{eqnarray*}
\end{exm}

The category $\C{cross}$ of crossed modules inherits a monoidal
structure $\otimes$ from the category of crossed complexes, since
it may be regarded as the full reflective subcategory of crossed
complexes concentrated in degrees one and two. The reflection
$\psi\colon\C{ccplx}\r\C{cross}$ sends a crossed complex
$$(C,\partial)=\left(\cdots\r C_n\st{\partial}\To C_{n-1}\r\cdots\r C_3 \st{\partial}\To C_2\st{\partial}\To C_1\right)$$
to the crossed module
\begin{equation}\label{lapsi}
\psi(C,\partial)=\left(\cdots\r 0\st{\partial}\To
0\r\cdots\r 0 \st{\partial}\To C_2/\partial(C_3)\st{\partial}\To
C_1\right).
\end{equation}
The unit of the reflection $(C,\partial)\r\psi(C,\partial)$ is the
identity in degree $1$, the natural projection
$C_2\twoheadrightarrow C_2/\partial(C_3)$ in degree $2$, and the
trivial map in higher degrees. Obviously any morphism from
$(C,\partial)$ to a crossed module factors uniquely through
$(C,\partial)\r\psi(C,\partial)$, so $\psi$ is indeed a
reflection.

The following lemma illustrates the rigidity of monoids in the category of crossed
modules of groups.

\begin{lem} {\rm (1)} Let $C$ be a crossed complex of groups and $\mu:C\otimes C\to C$
  a unital morphism. Then the induced morphism $\psi\mu:\psi C\otimes
  \psi C\to \psi C$ is a monoid structure.

{\rm (2)} Let $f:C\to C'$ be a morphism of crossed complexes of groups
  which preserves given unital morphisms $\mu:C\otimes C\to C$ and
  $\mu':C'\otimes C'\to C'$ up to some homotopy. Then $\psi f:\psi C\to\psi
  C'$ is a strict monoid homomorphism.
\end{lem}
\begin{proof}
(1) Since the only degree 0 element of $C$ is the unit, and $\mu$ is unital, the associativity relation
$\mu(\mu(a\otimes b)\otimes c)=\mu(a\otimes\mu(b\otimes
c))$ holds if the degree of $a$, $b$ or $c$ is 0.
If not, the total degree is at
least 3 and the relation is trivial on $\psi C$.

(2) Write $a_i$, $b_i$, $a_i'$, $b_i'$ for elements of $C_i$ and $C_i'$, $i\ge0$.
Since all the maps are unital,
$\mu'(fa_i\otimes fb_j)=f\mu(a_i\otimes b_j)$
if $i$ or $j=0$. It remains to show that this relation holds in the
crossed module $\psi
C'$ for $i=j=1$ also.

The homotopy will be given by a degree one function $h:C\otimes C\to
    C'$ satisfying a certain derivation formula
and an analogue of $\partial h+h\partial=\mu'(f\otimes f)-f\mu$, see e.g.~\cite{chII,tpcc}.

Clearly
    $\partial h(a_i\otimes b_j)=0$ for $\{i,j\}=\{0,1\}$, and
furthermore
the tensor
    product relations in $C'$ say that $\partial'(a_2'\otimes
    b_1')=\partial' a_2'\otimes b_1'-a_2'+{a_2'}^{b_1'}$.
In $\psi C'$ we can
    therefore deduce that $C_1'$ acts trivially
    on the elements  $a_2'=h(a_i\otimes b_j)$ for
$\{i,j\}=\{0,1\}$.
By the derivation property it now follows that in fact
    $h\partial(a_1\otimes b_1)=0$ in $\psi C'$, and so $\mu'(fa_1\otimes
    fb_1)=f\mu(a_1\otimes b_1)$ here also.
\end{proof}

\begin{cor}\label{salvo2}
{\rm (1)} Let $M\times M\to M$ be a strictly unital multiplication, where
$M$ is one of the following:
\begin{itemize}
\item a reduced simplicial set,
\item a reduced $CW$-complex,
\item a bisimplicial set with $M_{0,0}=\set{\text{point}}$.
\end{itemize}
Then  $\psi\pi M$, $\psi\pi_{CW}M$ or $\psi\Pi M$ respectively is a
monoid in the category of crossed modules.

{\rm (2)} Let $N\times N\to N$  be another such structure and $f\colon
M\r N$ a morphism which preserves multiplication up to a homotopy. Then
$f$ induces a strictly multiplicative homomorphism between the respective monoids in the category of crossed modules.
\end{cor}

Monoids in the category of crossed modules of groups are also termed reduced
$2$-modules, reduced $2$-crossed modules and strict braided categorical
groups, see \cite{2cm,2tils,3dnac,btc}. Commutative monoids are
similarly termed
stable crossed modules, stable $2$-modules and strict symmetric
categorical groups, see \cite{2cm,2tils,ccscg}.

We recall now the usual definition of these
concepts, following \cite{3dnac} and \cite{ccscg} up to a change of conventions.

\begin{defn}\label{larga}
A \emph{reduced $2$-module} is a crossed module $\partial\colon M\r N$ together with a map
$$\grupo{\cdot,\cdot}\colon N\times N\To M$$
satisfying the following identities for any $m,m'\in M$ and $n,n',n''\in N$.
\begin{enumerate}
\item $\partial\grupo{n,n'}=[n',n]$,
\item $m^n=m+\grupo{n,\partial(m)}$,
\item $\grupo{n,\partial(m)}+\grupo{\partial(m),n}=0$,
\item $\grupo{n,n'+n''}=\grupo{n,n'}^{n''}+\grupo{n,n''}$,
\item $\grupo{n+n',n''}=\grupo{n',n''}+\grupo{n,n''}^{n'}$.
\end{enumerate}
Moreover, $(\partial,\grupo{\cdot,\cdot})$ is a \emph{stable $2$-module} if (1), (2), (4) and
\begin{enumerate}\setcounter{enumi}{5}
\item $\grupo{n,n'}+\grupo{n',n}=0$
\end{enumerate}
are satisfied.

By (2), the action of $N$ on $M$ is completely determined by the
bracket $\grupo{\cdot,\cdot}$. The first crossed module axiom
(\ref{cm1}) is now redundant, and (\ref{cm2}) is equivalent
to
\begin{enumerate}\setcounter{enumi}{6}
\item $\grupo{\partial(m),\partial(m')}=[m',m]$,
\end{enumerate}
If we take (2) as a definition it is straightforward to check that
it does define a group action. Therefore we do not need to require
that $\partial$ is a crossed module, but just a homomorphism of
groups. Moreover, in the stable case (3) and (5) are redundant.
\end{defn}


\begin{lem}\label{reflex}
The category of stable quadratic modules is a full reflective
subcategory of the category of stable $2$-modules, given by those objects
$$C_0\times C_0\st{\grupo{\cdot,\cdot}}\To C_1\st{\partial}\To C_0$$
which satisfy
\begin{equation}\label{ley}
\grupo{c,[c',c'']}=0;\quad c, c', c''\in C_0.
\end{equation}
\end{lem}

\begin{proof}
This proof has two steps. We first identify the category of stable
quadratic modules with the full subcategory of $2$-modules
satisfying (\ref{ley}) and afterwards we show that this
subcategory is reflective.

We claim that a stable quadratic module $C_*$ yields a stable
$2$-module
\begin{equation*}\tag{a}
\xymatrix{C_0\times C_0\ar[rr]^-{(c,d)\mapsto
\ina{c}\otimes\ina{d}}\ar@/_20pt/[rrr]_{\grupo{\cdot,\cdot}}&&C^\abb_0\otimes
C^\abb_0\ar[r]^-\omega& C_1\ar[r]^{\partial}& C_0.}
\end{equation*}
This bracket satisfies (\ref{ley}) since
\begin{equation*}
\grupo{c,[c',c'']}\quad=\quad\omega(\ina{c}\otimes\ina{[c',c'']})\quad=\quad
\omega(\ina{c}\otimes0)\quad=\quad0.
\end{equation*}
Axioms (1), (6) and (7) in Definition \ref{larga} follow
immediately from Definition \ref{ob}, and axiom (4) is a
consequence of the following equations.
\begin{eqnarray*}
\grupo{c,c'+c''}&=&\omega(\ina{c}\otimes(\ina{c'}+\ina{c''}))\\
&=&\omega(\ina{c}\otimes\ina{c'})+\omega(\ina{c}\otimes\ina{c''})\\
&=&\grupo{c,c'}+\grupo{c,c''}\\
\text{{\small by (\ref{ley})}}\qquad&=&
\grupo{c,c'}+\grupo{c'',[c',c]}+\grupo{c,c''}\\
\text{{\small by Definition \ref{larga}
(1)}}\qquad&=&\grupo{c,c'}+\grupo{c'',\partial\grupo{c,c'}}+\grupo{c,c''}\\
\text{{\small by Definition \ref{larga}
(2)}}\qquad&=&\grupo{c,c'}^{c''}+\grupo{c,c''}.
\end{eqnarray*}
Therefore (a) is actually a stable $2$-module.

Conversely, let us see that a stable $2$-module
$$C_0\times C_0\st{\grupo{\cdot,\cdot}}\To C_1\st{\partial}\To C_0$$
satisfying (\ref{ley}) can be obtained from a stable quadratic
module as in (a). Indeed (\ref{ley}) and Definition \ref{larga}
(6) imply that $\grupo{\cdot,\cdot}$ factors through
$C_0^\abb\times C_0^\abb$. Moreover, by (\ref{ley}) and Definition
\ref{larga} (1) and (2) the elements of $C_0$ act trivially on the
image of $\grupo{\cdot,\cdot}$, therefore $\grupo{\cdot,\cdot}$ is
bilinear by  (4) and (5) in Definition \ref{larga}, so
$\grupo{\cdot,\cdot}$ factors through $C_0^\abb\otimes C_0^\abb$.
The factorization $\omega\colon C_0^\abb\otimes C_0^\abb\r C_1$ of
$\grupo{\cdot,\cdot}$ together with $\partial\colon C_1\r C_0$
define a stable quadratic module by Definition \ref{larga} (1),
(6) and (7), see Definition \ref{ob}.

We now prove that the subcategory of stable $2$-modules satisfying
(\ref{ley}) is reflective. Let
\begin{equation*}\tag{b}
N\times N\st{\grupo{\cdot,\cdot}}\To
M\st{\partial}\To N
\end{equation*}
be now an arbitrary stable $2$-module, let $N^\ni=N/[[N,N],N]$,
and let $P\subset M$ be the normal subgroup generated by the
subset $\grupo{N,[N,N]}\subset M$. We denote by $q\colon
N\twoheadrightarrow N^\ni$ and $q'\colon M\twoheadrightarrow M/P$
the natural projections. By Definition \ref{larga} (1) and (6) the
following diagram of solid arrows can be completed to a
commutative diagram in a unique way
\begin{equation*}\tag{c}
\xymatrix{N\times N\ar[r]^-{\grupo{\cdot,\cdot}}\ar@{->>}[d]_{q\times q}&M\ar[r]^-\partial\ar@{->>}[d]_{q'}&N\ar@{->>}[d]^q\\
N^\ni\times
N^\ni\ar@{-->}[r]_-{\grupo{\cdot,\cdot}}&M/P\ar@{-->}[r]_-\partial&N^\ni}
\end{equation*}
The subdiagram with dashed arrows is a stable $2$-module (the
axioms follow from the fact that (b) is a stable $2$-module and
that the vertical arrows in (c) are surjective), and this stable
$2$-module satisfies (\ref{ley}) by the way in which $P$ has been
defined. Moreover, it is immediate to check that any morphism from
(b) to a stable $2$-module factors through the projection (c) in a
unique way. This shows that the subcategory of stable $2$-modules
satisfying (\ref{ley}) is reflective and (c) is the unit of the
reflection.
\end{proof}

The reflection functor from stable $2$-modules to stable quadratic
modules constructed in the proof of Lemma \ref{reflex} will be
denoted by $\phi\colon\C{s2mod}\r\C{squad}$.
\medskip

Another nice feature of monoids in the category of
crossed modules of groups is that the property
of being commutative is preserved by weak equivalences.

\begin{lem}\label{escom}
Let $C\st{\sim}\r D$ be a morphism of reduced $2$-modules which is a weak equivalence. Then $C$ is
stable if and only if $D$ is.
\end{lem}

\begin{proof}
The operation
$\grupo{\cdot,\cdot}$ induces a natural quadratic function
$$\coker\partial\To\ker\partial\colon x\mapsto\grupo{x,x},$$
the $k$-invariant of $C$.
By using the properties of $\grupo{\cdot,\cdot}$ it is easy to see that $C$ is stable if and only if this quadratic function is indeed a group
homomorphism. Therefore the property of being stable is preserved under weak equivalences between reduced
$2$-modules.
\end{proof}

\begin{rem}\label{saleeq}
One
can obtain a stable $2$-module from an $(n-1)$-reduced
simplicial group $G$, $n\geq 2$, by using the following truncation
of the Moore complex $N_*G$
$$N_{n+1}G/d_0(N_{n+2}G)\st{d_0}\To N_nG=G_n.$$
The bracket is defined by
$$\grupo{x,y}=[s_1(x),s_0(y)]+[s_0(y),s_0(x)],\;\; x,y\in G_n.$$
This stable quadratic module will be denoted by
$\mu_{n+1}G$. If $G$ is only $0$-reduced this formula defines
a reduced $2$-module $\mu_2G$. Compare \cite{2cm,3dnac,ccscg}.

If $C\cong\mu_{n+1}G$ for an $(n-1)$-reduced free simplicial group
$G$, $n\geq 2$, then the natural morphism $C\r\phi C$ (i.e. the
unit of the reflection $\phi$ given by diagram (c) in the proof of
Lemma \ref{reflex}) is a weak equivalence. This is a consequence
of Curtis's connectivity result in \cite{srbhh} (which implies
that we can divide out weight three commutators in $G$ and still
obtain the same $\pi_n$ and $\pi_{n+1}$) since
$\mu_{n+1}(G/[[G,G],G])\cong\phi C$ and the natural morphism
$C\r\phi C$ is given by taking $\mu_{n+1}$ on the natural
projection $G\twoheadrightarrow G/[[G,G],G]$, compare
\cite[IV.B]{ch4c}.\ \ In order for $C$ to be such a truncation it
is enough that the lower-dimensional group of $C$ is free. Indeed
suppose that $E$ is the basis of the lower-dimensional group of
$C$. By \cite{2cm} there exists
an $(n-1)$-reduced simplicial group $G$
whose Moore complex is given by $C$ concentrated in dimensions $n$
and $n+1$.
In particular, $G_n=\grupo{E}$ is the free group with basis $E$.
By ``attaching cells'' one can construct a free resolution of $G$
(i.e.\ a cofibrant replacement) given by a weak equivalence
$G'\st{\sim}\r G$ in the category of simplicial groups which is
the identity in dimensions $\leq n$. Then
$\mu_{n+1}G'\cong\mu_{n+1}G\cong C $. As a consequence we observe
that the reflection $\phi$ preserves weak equivalences between
objects with a free low-dimensional group.
\end{rem}

Let $\ho\C{Spec}_0$ be the homotopy category of connective spectra of simplicial sets, and let
 $\ho\C{Spec}_0^1$ be the full subcategory of spectra with
trivial homotopy groups in dimensions other than $0$ and $1$.

\begin{lem}\label{equi}
There is a functor
$$\lambda_0\colon\ho\C{Spec}_0\To\ho\C{squad}$$
together with natural isomorphisms
$$\pi_i\lambda_0X\cong\pi_iX,\;\;i=0,1,$$ which induces an equivalence of
categories
$$\lambda_0\colon\ho\C{Spec}_0^1\st{\sim}\To\ho\C{squad}.$$
Moreover, for any connective spectrum the first Postnikov invariant of $X$ coincides with the $k$-invariant of
$\lambda_0X$.
\end{lem}

\begin{proof}
Stable quadratic modules, stable crossed modules and stable $2$-modules are known to be algebraic models of the $(n+1)$-type of an
$(n-1)$-reduced simplicial set $X$ for $n\geq 3$, see \cite{ch4c,2cm,2tils,ccscg}.
All these approaches are essentially equivalent, and they encode the first $k$-invariant as stated above. For
example, if $X$ is an $(n-1)$-reduced simplicial set, $n\geq 3$, then $\mu_nG(X)$
is such a model for the $(n+1)$-type of $X$. Here we use the Kan loop group $G(X)$. Its projection to stable
quadratic modules $\phi\mu_nG(X)$ is also a model for  the $(n+1)$-type of $X$ since $G(X)$ is free, see Remark \ref{saleeq} above.

The $1$-type of a connective spectrum $X$ of
simplicial sets is completely
determined by the $4$-type of the third simplicial set $Y_3$ of a fibrant replacement (in particular an
$\L$-spectrum) $Y$ of $X$. We can always assume that $Y_3$ is $2$-reduced. Therefore we
can define the functor $\lambda_0$ above as follows. Each spectrum $X$ is sent by $\lambda_0$ to
$\phi\mu_3G(Y_3)$.
\end{proof}

\begin{lem}\label{duro}
Given a $1$-reduced simplicial set $X$ there is a natural isomorphism
of monoids in crossed modules
$\psi\ul{\Omega}_{\C{Crs}}X\cong\mu_2G(X)$.
\end{lem}

\begin{proof}
Both $\psi\ul{\Omega}_{\C{Crs}}X$ and $\mu_2G(X)$ are models for
the $2$-type of the loop space of $X$, and moreover they have the
same low-dimensional group $\grupo{X_2-*}$, the free group with
basis $X_2-*$.

Using the presentation of $\ul{\Omega}_{\C{Crs}}X$ as a monoid in the
category of crossed complexes given in
\cite[Theorem 2.8]{tcc} and the convention followed by May
\cite[Definition 2.6.3]{saat} for the definition of $G(X)$, a
natural isomorphism
$\chi\colon\psi\ul{\Omega}_{\C{Crs}}X\cong\mu_2G(X)$
can be described on the monoid generators as follows. Given $x_2\in X_2$, let
$\chi(x_2)=x_2$, and given $x_3\in X_3$,
$$\chi(x_3)=-s_1d_2(x_3)+x_3-s_2d_3(x_3)+s_1d_3(x_3).$$
This is the identity in low-dimensional groups.
In order to check that it indeed defines an isomorphism in the upper groups one can use
the presentation of $\ul{\Omega}_{\C{Crs}}X$ in \cite{tcc},
and a computation
of
the Moore complex
of $G(X)$ in low dimensions by using the Reidemeister-Schreier method, see \cite[18]{cdhg} and \cite{mks}.
\end{proof}

In the statement of the following lemma we use the Moore loop complex
functor $\L$ on the category of fibrant simplicial sets.
Given a $1$-reduced Kan complex $X$, define $\L X$ by
$$(\L X)_{n}\;\;=\;\;\ker [d_{n+1}\colon
X_{n+1}\r X_n]$$
in the category of pointed sets; compare \cite[2.9]{sht}, \cite[Definition 23.3]{saat}.
The face and degeneracy operators are restrictions of the
operators in $X$. If $X$ is a simplicial group then so is $\L X$.

Recall that the natural simplicial map
$$\tau_X:\L X\r GX$$
given by $(\L X)_n\subset X_{n+1}\to\grupo{X_{n+1}}\twoheadrightarrow\grupo{X_{n+1}-s_0X_n}$ is a
homotopy equivalence when $X$ is a 1-reduced Kan complex.
The composite of $\pi_n\tau_X$ with
$\pi_{n+1}X\cong\pi_n\L X$ coincides with the connecting map
$\delta\colon\pi_{n+1}X\cong \pi_nGX$ in the path--loop
group fibration $GX\to
EX\to X$.

\begin{lem}\label{kk}
For any $2$-reduced Kan complex $X$ there is a natural weak
equivalence of simplicial groups $\sigma\colon G(\L X)\st{\sim}\r\L G(X)$.
\end{lem}

\begin{proof}
For all $n\geq 0$ we have
$$G_n(\L X)\cong \grupo{(\Omega X)_{n+1}-s_0(\Omega X)_n},$$
$$(\Omega G(X))_{n}\subset G_{n+1}(X)\cong \grupo{X_{n+2}-s_0X_{n+1}}.$$
The homomorphisms $\sigma_n:G_n(\L X)\r (\Omega G(X))_{n}$ are the unique possible homomorphisms compatible with the
inclusions $(\Omega X)_k\subset X_{k+1}$, $k\geq 0$, in the obvious way.
Since $\sigma\circ \tau_{\L X}=\L \tau_X:\L\L X\r\L GX$, the map
$\sigma$ is a weak equivalence.
\end{proof}

Now we are ready for the proof of the main theorem of this paper.

\begin{proof}[Proof of Theorem \ref{main}]
The coproduct in $\C{C}$ gives rise to a $\Gamma$-space $A$ in the
sense of Segal \cite{cct} with $A(\mathbf{1})=\abs{\diag\ner
wS.\C{C}}$, see \cite[Section 4, Corollary]{akttsI}. See also
the proof of Lemma \ref{step-by-step} below for further details on
$\ner wS.\C{C}$. The spectrum of topological spaces
$A(\mathbf{1})$, $BA(\mathbf{1})$, $B^2A(\mathbf{1})$,\dots
associated to $A$ is an $\Omega$-spectrum since $\diag\ner
wS.\C{C}$ is reduced. The $\Omega$-spectrum defining $K\C{C}$ is
obtained from the spectrum of $A$ by shifting the dimensions by
$+1$, i.e.\ $K\C{C}$ is given by
$$\Omega A(\mathbf{1}), A(\mathbf{1}),
BA(\mathbf{1}), B^2A(\mathbf{1}), \dots.$$

A particular choice of the coproduct $A\vee B$ of any pair of objects $A, B$ in
$\C{C}$ induces a product in $\ner wS.\C{C}$. We choose $A\vee *=A=*\vee A$ so that this product is strictly unital as in Corollary \ref{salvo2}.  The structure weak
equivalence
\begin{equation*}\tag{a}
\abs{\diag\ner
wS.\C{C}}\st{\sim}\To\Omega BA(\mathbf{1})
\end{equation*}
is a morphism of $H$-spaces up to homotopy.

We can replace $BA(\mathbf{1})$ and $B^2A(\mathbf{1})$ by homotopy equivalent spaces $\abs{Y_2}$, $\abs{Y_3}$
which are realizations of a $1$-reduced fibrant simplicial set $Y_2$ and a $2$-reduced fibrant simplicial set $Y_3$, respectively.
As a consequence we obtain a replacement for (a) consisting of a homotopy equivalence of pointed $CW$-complexes
\begin{equation*}\tag{b}
\abs{\diag\ner
wS.\C{C}}\st{\sim}\To\ul{\Omega}_{\C{FTop}} Y_2.
\end{equation*}
Here $\ul{\Omega}_{\C{FTop}} Y_2$ is the model for $\Omega\abs{Y_2}$ in \cite[Theorem 2.7]{tcc}. The $CW$-complex
$\ul{\Omega}_{\C{FTop}} Y_2$ is a monoid and the map (b) satisfies the hypotheses of Corollary \ref{salvo2}.

In order to define $\lambda_0K\C{C}$ as $\phi\mu_3G(Y_3)$ we choose an $\L$-spectrum $Y$ in the category of simplicial sets
representing $K\C{C}$ with $Y_2$
and $Y_3$ the simplicial sets chosen above.

Combining the results above we obtain the following weak equivalences of stable $2$-modules.

\[\begin{array}{rcll}
\psi\Pi\ner wS.\C{C}
&\st{\sim}\r&\psi\pi\diag\ner wS.\C{C}
&\text{(Theorems \ref{biEZ} and \ref{laxmon})}
\\
&=&\psi\pi_{CW}\abs{\diag\ner wS.\C{C}}
\!\!
&\text{ }
\\
&\st{\sim}\r&\psi\pi_{CW}\ul{\Omega}_{\C{FTop}} Y_2
&\text{(b)   }
\\
&\cong&\psi\ul{\Omega}_{\C{Crs}} Y_2
&\text{\cite[proof of Proposition 2.11]{tcc}   }
\\
&\cong&\mu_2 G(Y_2)
&\text{(Lemma \ref{duro})   }
\\
&\st{\sim}\r&\mu_2G(\Omega Y_3)
&\text{(Induced by $Y_2\st{\sim}\r\Omega Y_3$)    }
\\
&\st{\sim}\r&\mu_2\Omega G(Y_3)
&\text{(Lemma \ref{kk})   }
\\&=&\mu_3G(Y_3).
\end{array}\]
Here we use Lemma \ref{escom} to derive that not only $\mu_3G(Y_3)$ but all these reduced $2$-modules are indeed stable.

By Remark  \ref{saleeq}
we know that $\phi$ preserves weak equivalences between these stable $2$-modules, since they have free lower-dimensional group, so that
$$\phi\psi\Pi\ner
wS.\C{C}\st{\sim}\To\lambda_0 K\C{C}.$$
Now Theorem \ref{main} follows from Lemma \ref{step-by-step} below.
\end{proof}

\begin{lem}\label{step-by-step}
Let $\C C$ be a Waldhausen category and let $X=\ner wS.\C{C}$, the bisimplicial
set given by the nerve of the simplicial category $wS.\C C$. Then
there is an identification of stable quadratic modules
$$\D{D}_*\C{C}\;=\;\phi\psi\Pi X$$
which arises from an identification of the generators and relations on
both sides.
\end{lem}
\begin{proof}
We will give a presentation for the crossed module
$\psi\Pi(X)_2\st\partial\longrightarrow\psi\Pi(X)_1$,
and write down the commutative monoid structure
$$\langle\cdot,\cdot\rangle\colon\psi\Pi(X)\otimes\psi\Pi(X)\to \psi\Pi(X)$$
arising from the coproduct on $\C C$. This will give a presentation for
the associated stable quadratic module $\phi\psi\Pi X$ which
coincides with the presentation for our model $\D{D}_*\C{C}$, and we
will recover explicitly the relations (1)--(9) of Definition \ref{LA}.

Recall from \cite{akttsI} that Waldhausen's construction $wS.\C C$ is given by the simplicial object
$$\xymatrix@!=5em{
\{{*}\}=wS_0\C{C}\ar@/^4ex/[r]^{s_0}&\ar[l]<-0.4ex>_-{d_0}\ar[l]<0.4ex>^-{d_1}
wS_1\C{C}\ar@/^3ex/[r]<0.8ex>^{s_0}\ar@/^3ex/[r]<0.1ex>_{s_1}
&\ar[l]<-0.5ex>\ar[l]<0.2ex>\ar[l]<0.9ex>^{d_0,d_1,d_2}
wS_2\C{C}\ar@/^3ex/[r]^{s_0,s_1,s_2}&\ar[l]^{d_0,d_1,d_2,d_3}
wS_3\C{C}\ar@/^3ex/[r]&\ar[l]
\cdots}
$$
where each $wS_m\C C$ is a category whose objects $a$ are the sequences of cofibrations
$$
A_1\rightarrowtail A_2\rightarrowtail A_3\rightarrowtail\cdots\rightarrowtail A_{m-1}\rightarrowtail A_m
$$
and associated cofiber sequences
$A_j/A_i\rightarrowtail A_k/A_i\twoheadrightarrow A_k/A_j$
for $0\le i<j<k\le m$, with $A_0=*$.
For $i\neq0$ the (horizontal)
degeneracy operator $s_i$ inserts an identity map $A_i\r A_i$
and the face operator $d_i$ omits the object $A_i$,
while
$s_0$ inserts $*\rightarrowtail A_1$ and
$d_0$ replaces each $A_i$ by $A_i/A_1$, omitting $A_1$.

Morphisms $a\st\sim\r a'$ in $wS_m\C C$ are commutative diagrams which
are (levelwise) weak equivalences. Thus in the nerve $X$ we have sets of $(m,n)$-simplices,
$$X_{m,n}\;=\;\textrm{Ner}(wS_m\C{C})_n\,,$$
which consist of all composable strings of weak equivalences,
$$
x_{m,n}\;=\;\left(a\st\sim\longrightarrow a' \st\sim\longrightarrow a''
\st\sim\longrightarrow \cdots \st\sim\longrightarrow
a^{(n)}\right),
$$
between objects $a^{(i)}$ in $wS_m\C{C}$. We refer to \cite[1.4]{gj} for further details on the nerve
of a category.

Since $X$ is a {\em horizontally reduced} bisimplicial set, with $X_{0,n}\cong\{*\}$ for all $n$, Lemma
\ref{presen} will give a presentation of the total crossed
complex $\Pi(X)$.
In particular, the generators of $\Pi(X)$ in degree $k$
are the $(m,n)$-simplices of $X$ with total dimension $m+n=k$.
We need only consider degrees $k=1,2,3$,
since the crossed module $\psi\Pi(X)$ is the quotient of
$\Pi(X)_2\st{\partial_2}\longrightarrow\Pi(X)_1$
by the extra relation
\begin{equation}\label{d3zero}\partial_3(x)\;=\;0\end{equation}
in $\Pi(X)_2$, for each generator $x$ in $\Pi(X)_3$, see
(\ref{lapsi}). Note that degree $k=1$ will correspond to
$\D{D}_0\C{C}$, and $k=2$ to $\D{D}_1\C{C}$.

For $k=1$ we must consider the sets $X_{0,1}$ and $X_{1,0}$
but we can ignore the former since $X$ is horizontally reduced and the
first relation of Lemma \ref{presen} says $x_{0,1}=0$.
Elements $x_{1,0}\in X_{1,0}$ are
just objects $A$ in $\C C$, and so we can identify
generators in this degree with the symbols $[A]$.

For $k=2$ we consider the sets
$X_{1,1}$ and $X_{2,0}$, since again $X_{0,2}$ is degenerate.
Elements $x_{1,1}\in X_{1,1}$ are the weak equivalences between objects
in $\C C$, and elements  $x_{2,0}\in X_{2,0}$ are objects
of the category $wS_2\C C$ as above. Thus we can identify
generators in this degree with the symbols $[A\st\sim\r A']$ and
$[A\rightarrowtail B\twoheadrightarrow B/A]$ respectively.
\begin{figure}[h]\[\entrymodifiers={+0}
\begin{array}{c}\xymatrix@!R=3pc@!C=3pc{\ar@{-}[r]|*+{\scriptstyle A}&}\end{array}
\qquad\qquad
\begin{array}{c}
\xymatrix@!=2pc{\ar@{-}[r]|*+{\scriptstyle A'}="a"\ar@{-}[d]&\ar@{-}[d]\\
\ar@{-}[r]|*+{\scriptstyle A}="b"&\ar@{->}"b";"a"_{\sim}}
\end{array}
\qquad\quad
\begin{array}{c}
\xymatrix@!R=2.5pc@!C=2.5pc{&\\
\ar@{-}[r]|*+{\scriptstyle A}="a"\ar@{-}[ru]|*+{\scriptstyle B}="b"
&\ar@{-}[u]|*+{\scriptstyle B/A}="c"\ar@{>->}"a";"b"\ar@{->>}"b";"c"}
\end{array}
\]
\caption{Each 1-cell $x_{1,0}$ or 2-cell $x_{1,1}$ or $x_{2,0}$ corresponds
to a generator.}
\end{figure}
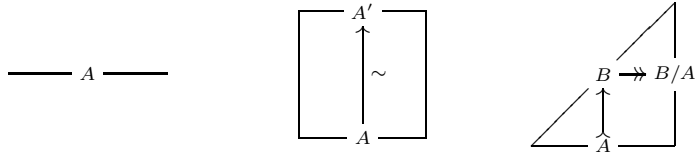

The second and third relations of Lemma~\ref{presen}
give the map $\partial$ on the degree 2 generators,
in terms of the horizontal and vertical simplicial face operators in $X$:
\begin{align*}
\tag{1}\partial([A\st\sim\r A'])
&={}-[d_0^v(A\st\sim\r A')]+[d_1^v(A\st\sim\r A')]\\
&={}-[A']+[A]\,, \\
\tag{2}\partial([A\rightarrowtail B\twoheadrightarrow B/A])
&={}-[d_1^h(A\rightarrowtail B\twoheadrightarrow B/A)]
    +[d_0^h(A\rightarrowtail B\twoheadrightarrow B/A)] \\
&\qquad\qquad\qquad\qquad\qquad\qquad\!\!{}
    +[d_2^h(A\rightarrowtail B\twoheadrightarrow B/A)] \\
&={}-[B]+[B/A]+[A]\,.
\end{align*}
Some of the above generators are trivial,
again by the first relation of Lemma \ref{presen}.
Explicitly, $x_{1,0}=0$ if it is $s^h_0$-degenerate,
$x_{1,1}=0$ if it is $s^v_0$-degenerate,
and $x_{2,0}=0$ if it is $s^h_0$- or $s^h_1$-degenerate:
\begin{align}
\tag{3} [*]&=[s_0^h(*)]\,=\,0\,,\\
\tag{4} [A\st{\!1_{\!A}\,}\r A]&=[s_0^v(A)]=0\,,\\
\tag{5} [{*}\rightarrowtail A\st{\!1_{\!A}\,}\twoheadrightarrow A]&=[s_0^h(A)]=0\,,\\
        [A\st{\!1_{\!A}}\rightarrowtail A\twoheadrightarrow{*}]&=[s_1^h(A)]=0\,.\nonumber
\end{align}

For $k=3$ we consider elements $x_{m,n}\in X_{m,n}$ for
$(m,n)=(1,2)$, $(2,1)$ and $(3,0)$.
An element $x_{1,2}$ is a composable pair of weak equivalences in $\C C$,
an element $x_{2,1}$ is a weak equivalence in $wS_2\C C$,
and an element $x_{3,0}$ is an object of $wS_3\C C$.
\begin{figure}[h]
\[\entrymodifiers={+0}\xymatrix@!R=2pc@!C=0.3pc{
 \arr@{-}[rrr]|*+{\scriptscriptstyle C}="c"&&&\\
&\arr@{.}[lu]&&&\arr@{-}[lu]\\
\arr@{-}[uu]\arr@{.}[ru]\arr@{-}[rrr]|*+{\scriptscriptstyle A}="a"&&&\arr@{-}[ru]\arr@{-}[uu]
\arr@{.}"2,2";"2,5"|!{"c";"a"}\hole|*+{\scriptscriptstyle B}="b"|!{"1,4";"3,4"}\hole
\arr@*{}@{.>}"a";"b"_*{\scriptscriptstyle\!\!\!\sim}
\arr@*{}@{.>}"b";"c"_*{\scriptscriptstyle\!\!\!\sim}
\arr@*{}@{->}"a";"c"^(0.3)*{\scriptscriptstyle\sim\!}
}
\!\!\!\! x_{1,2}\qquad\qquad\qquad\qquad\qquad x_{2,1}\qquad
\xymatrix@!C=1pc@!R=0.7pc{
\arr@{-}[rrr]|*+{\scriptscriptstyle B'}="bb"
\arr@{-}[rd]|*+{\scriptscriptstyle A'}="aa"&&&\\&
\arr@{-}[urr]|*+{\scriptscriptstyle B'/A'}="bbaa"\\\arr@{-}[uu]
\arr@{.}[rrr]|(0.1666666){\hole}|(0.3333333){\hole}
         |*+{\scriptscriptstyle B}="b"|(0.6666666){\hole}
\arr@{-}[rd]|*+{\scriptscriptstyle A}="a"&&&\arr@{-}[uu]\\&\arr@{-}[uu]
\arr@{-}[urr]|*+{\scriptscriptstyle B/A}="ba"
\arr@*{}@{>->} "aa";"bb"
\arr@*{}@{>.>}"a";"b"|{\hole}
\arr@*{}@{->>}  "bb";"bbaa"
\arr@*{}@{.>>} "b";"ba"
\arr@*{}@{->}  "a";"aa"^*{\scriptscriptstyle\sim}
\arr@*{}@{->}  "ba";"bbaa"^*{\scriptscriptstyle\sim}
\arr@*{}@{.>} "b";"bb"^(0.3)*{\scriptscriptstyle\sim}
|!{"2,2";"bbaa"}{\hole}}
\]
\vspace*{-23ex}

\[
\def\tetra#1#2#3#4#5{{\entrymodifiers={+0}\xymatrix@!R=3.5pc@!C=0.8pc{
#5\\#5\\#5\\ x_{3,0} #5\\#5
    \arr@{-}"#1";"#2"|*+{\scriptscriptstyle A}="a"
    \arr@{-}"#1";"#4"|*+{\scriptscriptstyle C}="c"
    \arr@{-}"#2";"#3"|*+{\scriptscriptstyle B/A}="ba"
    \arr@{-}"#2";"#4"|*+{\scriptscriptstyle C/A}="ca"
    \arr@{-}"#3";"#4"|*+{\scriptscriptstyle C/B}="cb"
    \arr@{.}"#1";"#3"|!{"a";"c"}\hole|*+{\scriptscriptstyle B}="b"|!{"#2";"#4"}\hole|!{"ba";"ca"}\hole
    \arr@{.>>}@*{}"b";"ba"|!{"#2";"#4"}\hole
    \arr@{.>>}@*{}"c";"cb"|!{"#2";"#4"}\hole
    \arr@{>.>}@*{}"a";"b"
    \arr@{>->} @*{}"a";"c"
    \arr@{>.>}@*{}"b";"c"
    \arr@{>->} @*{}"ba";"ca"
    \arr@{->>} @*{}"c";"ca"
    \arr@{->>} @*{}"ca";"cb"}}}
\tetra{3,1}{4,4}{3,6}{1,4}{&&&&&}\qquad\qquad
\]
\vspace*{-8ex}

\caption{Each 3-cell in $X$ gives a relation $\partial_3(x)=0$
  (equation \ref{d3zero}).}
\end{figure}
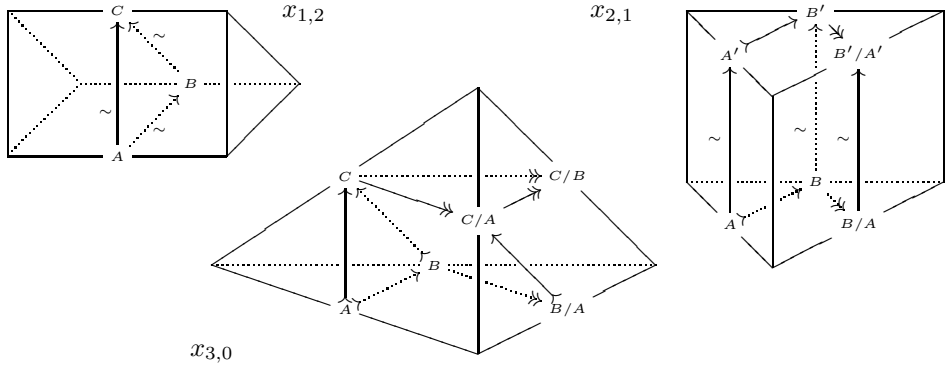

\noindent%
Thus by the fourth, fifth and sixth relations of Lemma \ref{presen} we
have
\begin{align*}
  \tag{6} 0&=\partial_3(x_{1,2}) ={}
-d_2^vx_{1,2}-d_0^vx_{1,2}+d_1^vx_{1,2}
\\&= {}-[A\st{\sim}\r B]-[B\st{\sim}\r C]+[A\st{\sim}\r C]\,,
\\
\tag{7} 0&=\partial_3(x_{2,1})  =
 d_2^hx_{2,1}
+{d_0^hx_{2,1}}^{d_2^hd_1^vx_{2,1}}
-d_1^vx_{2,1}
-d_1^hx_{2,1}
+d_0^vx_{2,1}
\\&={}
[A\st{\sim}\r A']+[B/A\st{\sim}\r B'/A']+\grupo{[A],-[B'/A']+[B/A]}
\\&
{}
\quad-[A \rightarrowtail B \twoheadrightarrow B/A]
-[B\st{\sim}\r B']
+[A'\rightarrowtail B'\twoheadrightarrow B'/A']\,,
\\
\tag{8} 0&=\partial_3(x_{3,0}) =
 d^h_2x_{3,0}
+{d^h_0x_{3,0}}^{(d_2^h)^2x_{3,0}}
-d^h_3x_{3,0}
-d^h_1x_{3,0}
\\&={}
  [A\rightarrowtail C\twoheadrightarrow C/A]
 +[B/A\rightarrowtail C/A\twoheadrightarrow C/B]
\\&\quad +\grupo{[A],-[C/A]+[C/B]+[B/A]}
 -[A\rightarrowtail B\twoheadrightarrow B/A]
 -[B\rightarrowtail C\twoheadrightarrow C/B]
\end{align*}
using Definition \ref{larga} (2) to rewrite the actions.

Finally, recall that by Corollary \ref{salvo2} the operation $\vee\colon X\times X\to X$ given
by the coproduct in $\C C$ induces a monoid structure on $\psi\Pi(X)$,
$$\langle\cdot,\cdot\rangle\colon
\psi\Pi(X)\otimes\psi\Pi(X)\st{b''}\longrightarrow
\psi\Pi(X\times X)\st{\vee}\longrightarrow
\psi\Pi(X)\,.$$
By the formula for the shuffle map $b''$ in Example \ref{ll} we have
\begin{align*}
\tag{9}&\langle[A],[B]\rangle={}-[s_0^h(A)\vee s_1^h(B)]+[s_1^h(A)\vee s_0^h(B)]
\\&\qquad=-[(*\rightarrowtail A\st1\r A)\vee (B\st1\r B\twoheadrightarrow *)]+
[(A\st1\r A\twoheadrightarrow *)\vee(*\rightarrowtail B\st1\r B)]
\\&\qquad=-[B\st{i_2}\rightarrowtail A\vee B\st{p_1}\twoheadrightarrow A]
+[A\st{i_1}\rightarrowtail A\vee B\st{p_2}\twoheadrightarrow B].
\end{align*}
Here $i_j$ and $p_j$ denote respectively the natural inclusions
and the projections of the two factors of the coproduct $A\vee B$,
$j=1,2$.
\end{proof}

\appendix

\section{Free stable quadratic modules and presentations}\label{append}

Let
$$U\colon\C{squad}\To\C{Set}\times\C{Set}$$
be the forgetful functor from stable quadratic modules to pairs of
sets defined by $U(C_*)=(C_0,C_1)$. The functor $U$ has a left
adjoint $F$, and a stable quadratic module $F(E_0,E_1)$ is called
a {\em free stable quadratic module} on the sets $E_0$ and $E_1$.
In order to give an explicit description of $F(E_0,E_1)$ we fix
some notation.

Given a set $E$ we denote by $\grupo{E}$ the \emph{free group}
with basis $E$, and by $\grupo{E}^{\abb}$ the \emph{free abelian
group} with basis $E$. The \emph{free group of nilpotency class
$2$} with basis $E$, denoted by $\grupo{E}^{\ni}$, is the quotient
of $\grupo{E}$ by triple commutators.

For any abelian group $A$ let $\hat{\otimes}^2A$ be the quotient
of the tensor square $A\otimes A$ by the relations $a\otimes
b+b\otimes a=0$, $a,b\in A$. The projection of $a\otimes b\in
A\otimes A$ to $\hat{\otimes}^2A$ is denoted by $a\hat{\otimes}b$.

Given a pair of sets $E_0$ and $E_1$, we write $E_0\cup
\partial E_1$ for the set whose elements are the symbols $e_0$ and
$\partial e_1$ for each $e_0\in E_0$, $e_1\in E_1$.

To define $F(E_0,E_1)$, consider the groups
\begin{eqnarray*}
F(E_0,E_1)_0&=&\grupo{E_0\cup \partial E_1}^\ni,\\
F(E_0,E_1)_1&=&\hat\otimes^2\grupo{E_0}^{\abb}\times
\grupo{E_0\times E_1}^\abb\times\grupo{E_1}^\ni.
\end{eqnarray*}
The structure homomorphisms of $F(E_0,E_1)$ are defined as
follows. Given $e_i,e_i',e_i''\in E_i$,
\begin{eqnarray*}
\partial(e_0\hat{\otimes}e_0',(e_0'',e_1),e_1')&=&[e_0',e_0]+[\partial
e_1,e_0'']+\partial e_1',\\
\grupo{e_0,e_0'}&=&(e_0\hat{\otimes}e_0',0,0),\\
\grupo{e_0,\partial e_1}&=&(0,e_0\otimes e_1,0)\;=\;-\grupo{\partial e_1,e_0},\\
\grupo{\partial e_1,\partial e_1'}&=&(0,0,[e_1',e_1]).
\end{eqnarray*}
In the language of \cite[IV.C]{ch4c} $F(E_0,E_1)$ is the
\emph{totally free} stable quadratic module with basis given by
the function $E_1\r\grupo{E_0\cup \partial E_1}^{nil}\colon
e_1\mapsto
\partial e_1$. Therefore $F(E_0,E_1)$ is indeed a stable quadratic module and, moreover, it
satisfies the required universal property, i.e. given a stable
quadratic module $C_*$ any pair of maps $E_i\r C_i$, $i=0,1$, can
be uniquely extended to a morphism $F(E_0,E_1)\r C_*$ in
$\C{squad}$.


It is now straightforward to define explicitly the stable
quadratic module $C_\ast$ presented by {\em generators} $E_i$ and
{\em relations} $R_i\subset F(E_0,E_1)_i$ in degrees $i=0,1$, by
\begin{eqnarray*}
C_0&=&F(E_0,E_1)_0/N_0,\\
C_1&=&F(E_0,E_1)_1/N_1.
\end{eqnarray*}
Here $N_0\subset F(E_0,E_1)_0$ is the normal subgroup generated by
the elements of $R_0$ and $\partial (R_1)$, and $N_1\subset
F(E_0,E_1)_1$ is the normal subgroup generated by the elements of
$R_1$ and $\grupo{F(E_0,E_1)_0,N_0}$. The structure homomorphisms
of $F(E_0,E_1)$ induce a stable quadratic module structure on
$C_*$ which satisfies the following universal property: given a
stable quadratic module $C_*'$, any pair of functions $E_i\r C_i'$
$(i=0,1)$ such that the induced morphism $F(E_0,E_1)\r C_*'$
annihilates $R_0$ and $R_1$ induces a unique morphism $C_*\r
C_*'$.

In \cite[IV.C]{ch4c} Baues also considers the \emph{totally free}
stable quadratic module $C_*$ with basis given by a function
$g\colon E_1\r\grupo{E_0}^{nil}$. In the language of this paper
$C_*$ is the stable quadratic module with generators $E_i$ in
degree $i=0,1$ and degree $0$ relations $\partial(e_1)=g(e_1)$ for
all $e_1\in E_1$.

\bibliographystyle{amsalpha}
\bibliography{Fernando}
\end{document}